\documentstyle{amltd2004}
\begin{document}
\annalsline{155}{2002}
\startingpage{553}
\def\bye{\end{document}}
 \font\tenrm=cmr10

\input amssym.def
\input amssym.tex
\def\lrp#1{\left(#1\right)}

\catcode`\@=11
\font\twelvemsb=msbm10 scaled 1100
\font\tenmsb=msbm10
\font\ninemsb=msbm10 scaled 800
\newfam\msbfam
\textfont\msbfam=\twelvemsb  \scriptfont\msbfam=\ninemsb
  \scriptscriptfont\msbfam=\ninemsb
\def\msb@{\hexnumber@\msbfam}
\def\Bbb{\relax\ifmmode\let\next\Bbb@\else
 \def\next{\errmessage{Use \string\Bbb\space only in math
mode}}\fi\next}
\def\Bbb@#1{{\Bbb@@{#1}}}
\def\Bbb@@#1{\fam\msbfam#1}
\catcode`\@=12

 \catcode`\@=11
\font\twelveeuf=eufm10 scaled 1100
\font\teneuf=eufm10
\font\nineeuf=eufm7 scaled 1100
\newfam\euffam
\textfont\euffam=\twelveeuf  \scriptfont\euffam=\teneuf
  \scriptscriptfont\euffam=\nineeuf
\def\euf@{\hexnumber@\euffam}
\def\frak{\relax\ifmmode\let\next\frak@\else
 \def\next{\errmessage{Use \string\frak\space only in math
mode}}\fi\next}
\def\frak@#1{{\frak@@{#1}}}
\def\frak@@#1{\fam\euffam#1}
\catcode`\@=12


\newcommand{\twosum}[2]{\sum_{\begin{array}{c} {\scriptstyle #1}\\
{\scriptstyle #2} \end{array}}}
\newcommand{\Q}{{\Bbb Q}}
\newcommand{\F}{{\Bbb F}}
\newcommand{\R}{{\Bbb R}}
\newcommand{\C}{{\Bbb C}}
\newcommand{\N}{{\Bbb N}}
\newcommand{\Z}{{\Bbb Z}}
\newcommand{\Qb}{\overline{\Q}}
\newcommand{\Proj}{{\Bbb P}}
\newcommand{\Grass}{{\Bbb G}}
\newcommand{\ep}{\varepsilon}
\newcommand{\lb}[1]{\mbox{\scriptsize\bf #1}}
\newcommand{\mod}[1]{\hspace{-2.9mm}\pmod{#1}}
\renewcommand{\b}[1]{{\bf #1}}
\newcommand{\cl}[1]{{\cal #1}}
\newcommand{\x}{{\bf x}}
\newcommand{\sfLambda}{\mathbold{\Lambda}}
\newcommand{\sfMu}{\mathbold{\Mu}}
\newcommand{\rank}{{\rm rank}}
\newcommand{\sLambda}{\mathbold{\Lambda}}

\title{The density of rational points\\ on curves and surfaces}
\shorttitle{Rational points on curves and surfaces} 
\author{D. R. Heath-Brown}


\section{Introduction}

Let $n\geq 3$ be an integer and let 
$F(\x)=F(x_{1},\ldots,x_{n})\in\Z[x_{1},\ldots,x_{n}]$ be an
absolutely irreducible form of degree $d$, producing a hypersurface
of dimension $n-2$ in $\Proj^{n-1}$.  This paper is primarily concerned
with the number of rational points on this hypersurface, of height at most
$B$, say.  In order to describe such points we choose representatives
$\x=(x_{1},\ldots,x_{n})\in\Z^{n}$ with the $x_{i}$ not all $0$, and
such that $\gcd(x_{1},\ldots,x_{n})=1.$  Moreover we shall insist that
if $i$ is the smallest index for which $x_{i}\not=0$, then $x_{i}>0$.
We shall define $Z_{n}$ to be the set of all such representatives
$\x$.  Our primary interest is then with the quantity 
$$N(B)=N(F;B)=\#\{\x\in Z_{n}:
\,F(\x)=0,\;\max_{1\leq i\leq n}|x_{i}|\leq B\}.$$

We begin with a rather trivial result.

\specialnumber{1}\proclaim{Theorem}
For any $n\geq 2$ we have
\begin{equation}
N(F;B)\ll B^{n-1}.
\end{equation}
This remains true if $F$ is allowed to have coefficients in $\Qb${\rm .}
\endproclaim

Here, and throughout the paper, the implied constant may depend on $n$
and $d$.  However where there is a dependence on $F$ we shall say so
explicitly.  The result shows
in particular that there is an
integer vector $\x$ with $F(\x)\not=0$ satisfying $|\x|\ll_{n,d}1$,
and this is a fact that we shall use repeatedly.  We do not claim that
Theorem 1 is new. 

It is trivial that the exponent $n-1$ above is best
possible, in the case $d=1$.  However for $d=2$ we have
\begin{equation}
N(B)\ll_{F,\ep} B^{n-2+\ep},
\end{equation}
for any $\ep>0$.  This can be proved by the circle method, for example.
Our next result is a version of this which is independent of $F$.  

\specialnumber{2}\proclaim{Theorem}
Let $F(\x)$ be a quadratic form of rank at least $3${\rm ,} in $n$ variables{\rm .}
Then
$$N(B)\ll_{\ep} B^{n-2+\ep},$$
for any fixed $\ep>0${\rm .}
\endproclaim

As with Theorem 1, this estimate is almost trivial.  Again we do not
claim that the result is new.

We will be interested in the extent to which one can
prove results of this kind when $d\geq 3$.  
Let us first consider the case $n=3$, corresponding to curves in
$\Proj^{2}$.  When $d\geq 3$ and the curve has genus $1$, we have
N\'{e}ron's result
\begin{equation}
N(B)\sim c_{F}(\log B)^{r/2}\ll_{F,\ep}B^{\ep},
\end{equation}
for any $\ep>0$, where $c_{F}$ is a positive constant depending on $F$,
and $r$ is the rank of the Jacobian of the curve.
For genus $2$ or more we even have
\begin{equation}
N(B)\ll_{F}1,
\end{equation}
by the celebrated theorem of Faltings [8].  Unfortunately it is hard to
produce versions of these results with a good explicit dependence on
$F$.  Nonetheless it has been shown by Pila [27], via quite different
methods, that
\begin{equation}
N(B)\ll_{\ep} B^{1+1/d+\ep},
\end{equation}
for $n=3$ and any $\ep>0$.  Indeed Pila shows in general that
\begin{equation}
N(B)\ll_{\ep} B^{n-2+1/d+\ep}.
\end{equation}
It is remarkable that these results are completely
independent of $F$.  Pila's estimates are deduced from a bound relating
to integral points on affine curves due to Bombieri and Pila [2].
(See also Pila [28].)
Bombieri and Pila showed that if 
$f(x,y)\in\Z[x,y]$ is an absolutely irreducible
polynomial of degree $d$, then
\begin{equation}
\#\{(x,y)\in\Z^{2}:\, f(x,y)=0,\, |x|,|y|\leq B\}\ll_{\ep}
B^{1/d+\ep}.
\end{equation}
Our principal strategy in this paper will be to generalize this latter
result.  In particular we shall consider a projective version of it,
and we shall replace the cube of side $2B$ by a more general box.
This will prove very convenient for applications.  We therefore
take $\b{B}=(B_{1},\ldots,B_{n})$ with each $B_{i}\geq 1$, and define
a counting function
$$N(\b{B})=N(F;\b{B})=\#\{\x\in Z_{n}:\,F(\x)=0,\,
|x_{i}|\leq B_{i},\,(1\leq i\leq 1)\}.$$
It will be convenient to write
$$V=\prod_{i=1}^{n}B_{i}$$
and
$$T=\max\left\{\prod_{i=1}^{n}B_{i}^{f_{i}}\right\},$$
with the maximum taken over all integer $n$-tuples
$(f_{1},\ldots,f_{n})$ for which the corresponding monomial
$$x_{1}^{f_{1}}\ldots x_{n}^{f_{n}}$$
occurs in $F(\b{x})$ with nonzero coefficient.  When $d\geq n\geq 3$ and
$F$ is nonsingular, we may bound $T$ from below as follows.  We 
order the variables $X_{i}$ so that $B_{1}\geq B_{2}\geq\ldots\geq
B_{n},$ and observe 
that some monomial $x_{1}^{d-1}x_{i}$ must occur in $F(\b{x})$, where
$1\leq i\leq n$.  Thus $T\geq B_{1}^{d-1}B_{n}\geq V^{d/n}$.  Indeed,
for a nonsingular ternary quadratic form, the same estimate $T\geq
V^{2/3}$ still holds.  To see this, observe as above that $T\geq
B_{1}B_{2}\geq V^{2/3}$ if there is a term in $x_{1}^{2}$ or
$x_{1}x_{2}$.  If neither of these is present there must be terms in 
both $x_{1}x_{3}$ and $x_{2}^{2}$, since $F$ is nonsingular.  In this
case we have $T\geq\max(B_{1}B_{3},B_{2}^{2})\geq V^{2/3}$.

Our principal result for curves is the following.

\specialnumber{3}\proclaim{Theorem}
Let $n=3$ and $\ep>0${\rm .}  If $F$ is irreducible over $\Q${\rm ,} then
\begin{equation}
N(F;B_{1},B_{2},B_{3})\ll_{\ep} T^{-d^{-2}}V^{d^{-1}+\ep}.
\end{equation}
In particular we have
\begin{equation}
N(F;B)\ll_{\ep} B^{2/d+\ep}.
\end{equation}
Moreover if $F$ is nonsingular we have
\begin{equation}
N(F;B_{1},B_{2},B_{3})\ll_{\ep} V^{2/(3d)+\ep}.
\end{equation}
\endproclaim

As with the result of Bombieri and Pila, we have estimates that are
completely independent of $F$.  In fact this arises through an
application of the following result, in which we write $||F||$ for 
the height of the form $F$, defined
as the maximum modulus of the coefficients of $F$.

\specialnumber{4}\proclaim{Theorem} \hskip-8pt
Let $F(\x)\in\Z[\x]$ be a form in $n=3$ variables{\rm ,} of degree~$d${\rm .}
Suppose that $F$ is irreducible over $\Q${\rm ,} and that the coefficients 
of $F$ are coprime{\rm .}  Then
either $N(F;B)\leq d^{2}$ or $||F||\ll B^{d(d+1)(d+2)/2}${\rm .}
\endproclaim

This enables us to absorb a dependence of the type $||F||^{\ep}$ in
the estimate, into the term $V^{\ep}$ (or $B^{\ep}$).  A similar 
technique can be applied to higher dimensional varieties; see  Sections 5, 6 and  8.

One should note that
the exponent in (1.9) is appreciably smaller than that in (1.5).
Moreover, if we take $B_{1}=B_{2}=B$ and $B_{3}=1$ in (1.8) we recover
the exponent $1/d$ of (1.7).  We may also observe that if
$F(x_{1},x_{2},x_{3})=x_{1}^{d}-x_{2}^{d-1}x_{3}$, then the solutions
$(m^{d-1}n,m^{d},n^{d})$ show that 
$$N(F;B)\gg B^{2/d},$$
so that (1.9) is, in a suitable sense, best possible.  Finally it should
be pointed out that we do
not require $F$ to be absolutely irreducible for Theorem 3.  Indeed
for forms which are irreducible over $\Q$ but reducible over
$\overline{\Q}$ a stronger estimate is a consequence of the following result.

\specialnumber{1}\proclaim{{C}orollary}
Theorem $3$ holds for any $F(\x)\in\overline{\Q}[x_{1},x_{2},x_{3}]$
which is irreducible over $\overline{\Q}$ and of degree $d${\rm .}
Indeed if $F$ is not a multiple of a rational form then $N(F;B)\leq d^{2}${\rm . }
\endproclaim

The first statement clearly follows from the second.  To prove the
latter one merely writes $F$ as a linear combination
$\sum \lambda_{i}F_{i}$ of rational forms $F_{i}$, with linearly
independent $\lambda_{i}$.  Some $F_{i}$ is not a
multiple of $F$, but all rational zeros of $F$ must satisfy
$F=F_{i}=0$. The result then follows by B\'{e}zout's theorem.

At this point we remark that we shall use the term `absolutely
irreducible' to describe a polynomial, or an equation, which is
irreducible over $\overline{\Q}$.  When we only say `irreducible', the
relevant field must be understood from the context.  When the relevant
field is $\overline{\Q}$ we shall use the two terms interchangeably.
However, in the
context of curves and higher dimensional varieties, we shall use the
phrase `irreducible' to mean irreducible over $\overline{\Q}$.

When the $B_{i}$ are unequal, Theorem 3 is new even in the case
$d=2$.  In the author's work [13; p.\  24], the estimate
$$N(F;B_{1},B_{2},B_{3})\ll_{\ep} V^{1/2},$$
given by [13; Lemma 2] was employed.  By substituting the bound (1.8)
we can strengthen [13; Theorem 2] as follows:

\specialnumber{2}\proclaim{{C}orollary}\hskip-8pt
Let $q$ be an integral ternary quadratic form with matrix~$M.$ Let
$\Delta=|\det(M)|,$ and assume that $\Delta\not=0.$  Write
$\Delta_{0}$ for the highest common factor of the $2\times 2$ minors of
$M.$  Then the number of primitive integer solutions of
$q(\b{x})=0$
in the box $|x_{i}|\leq R_{i}$ is 
$$\ll_{\ep} \left\{1+\lrp{\frac{R_{1}R_{2}R_{3}\Delta_{0}^{2}}{\Delta}}^{1/3+\ep}\right\}
d_{3}(\Delta)$$
for any $\ep>0${\rm .}
\endproclaim

In the original version the exponent $1/3+\ep$ was replaced by $1/2$.
We may of course replace $d_{3}(\Delta)$ by $(R_{1}R_{2}R_{3})^{\ep}$
if we wish, by virtue of Theorem 4.

One can also estimate the number of points on a curve in $\Proj^{3}$.

\specialnumber{5}\proclaim{Theorem}
Let $C$ be an irreducible curve in $\Proj^{3}${\rm ,} of degree $d${\rm ,} not
necessarily defined over the rationals{\rm .}
Then $C$ has
$O_{\ep}(B^{2/d+\ep})$ points $\x\in Z_{4}$ in the cube
$\max |x_i|\leq B$.
\endproclaim \pagebreak

This will be established by projecting $C$ onto a suitable plane, and
counting the points on the resulting plane curve.  If $C$ is nonplanar
one would expect to loose information by such a process.  However in our
applications of Theorem~5 we are usually unable to tell whether or not
$C$ is planar.

It is interesting to compare the estimates given by Theorems 3 and 5 with
those obtained very recently by Elkies [6].  Elkies' emphasis is on
algorithms for finding rational points.  Thus he shows in [6; Theorem 3]
that one can find the rational points of height at most $B$, on a
curve $C$ of degree $d$, in time $O_{C,\ep}(B^{2/d+\ep})$.  It follows
in particular that there are $O_{C,\ep}(B^{2/d+\ep})$ points to be found.
Elkies does not consider issues
of uniformity with respect to the curve, although it seems quite
plausible that his methods will yield a good dependence on the height
of $C$, or even complete independence as in the present work.  At
first sight the approach taken in the two papers is rather different,
but closer inspection reveals interesting parallels.  Indeed Elkies
goes on to examine the situation for varieties of higher dimension,
presenting a heuristic argument that produces the same exponents
$3/\sqrt{d}$ and $(n-1)d^{-1/(n-2)}$ which arise from Theorem 14 below.

We now discuss the case $n=4$, corresponding to surfaces in $\Proj^{3}$.
The example $F(\x)=x_{1}^{d}+x_{2}^{d}-x_{3}^{d}-x_{4}^{d}$, for which
all vectors $(a,b,a,b)$ are solutions, shows
that we may have
$N(B)\gg B^{2}$
even when $F$ is nonsingular.  It is thus natural to exclude trivial
solutions by defining $N_{1}(B)$ to count the same rational points as
does $N(B)$, but excluding any that lie on lines in the
surface $F(\x)=0$.  We may then conjecture that
\begin{equation}
N_{1}(B)\ll_{F,\ep}B^{1+\ep}
\end{equation}
for any $\ep>0$, as soon as $d\geq 3$.  In so far as the weaker bound (1.2)
has not hitherto been established for general forms, even in the cubic 
case, the above conjecture is a long way off.  We may observe that if
$d\geq 3$, the
surface
$$x_{1}^{d}+x_{2}^{d}-x_{2}^{d-2}x_{3}x_{4}=0$$
is absolutely irreducible, and contains no lines other than those in
the planes $x_{2}=0,\,x_{3}=0$ and $x_{4}=0$.  
However there are rational points
$(0,ab,a^{2},b^{2})$, which show that $N_{1}(B)\gg B$ in this case.
Thus the exponent $1$ in (1.11) would be best possible.

We shall make some modest progress towards the above conjecture by
establishing the following result.

\specialnumber{6}\proclaim{Theorem}
For any  absolutely irreducible form
$F(\x)\in\overline{\Q}[x_{1},\ldots,x_{4}]$ of degree $3$ or more{\rm ,} 
we have 
$$N_{1}(F;B)\ll_{\ep}B^{52/27+\ep}.$$
\endproclaim

An inspection of the proof shows that the exponent $52/27$ may be
replaced by $17/9$ when $F$ has degree $4$ or more.  However 
we can improve substantially on this for large values of $d$, as follows.

\specialnumber{7}\proclaim{Theorem}
For any absolutely irreducible form
$F(\x)\in\overline{\Q}[x_{1},\ldots,x_{4}]$
of degree $d${\rm ,}
we have 
$$N_{1}(F;B)\ll_{\ep}B^{1+3/\sqrt{d}+\ep}.$$
\endproclaim

Theorem 6 answers questions raised by the author [14], by showing that
points on any lines in the surface $F=0$ that are defined over $\Q$ will
dominate $N(B)$.  Surfaces of the type
$G(x_{1},x_{2})=G(x_{3},x_{4})$, where $G$ is a binary form, have
been investigated fairly extensively.  Thus Hooley [16], [22] has
shown, in effect, that $N_{1}(B)=o(B^{2})$ 
when $G$ is a cubic form, and also [19] when
$G$ is a quartic form of the special type $G=ax^{4}+bx^{2}y^{2}+cy^{4}$.  
For binary forms of degree $d\geq 5$, the most general case that has been
covered is that of forms of the type $G=Ax^{d}+By^{d}$, which have
been handled by Bennett, Dummigan and Wooley [1].  There has however
been much work on the forms $G=x^{d}+y^{d}$, to which we shall allude later.
The sieve methods used by Hooley [16], [22] save a power of $\log B$
relative to $B^{2}$, whereas the other techniques used hitherto, which
trace their origins to Hooley's work [17] on sums of 4 cubes, save a
power of $B$.  

As a consequence of Theorem 6, we can show, in the spirit of the
above works, that most numbers represented by a binary form $G$ have
essentially only one representation.  To make this precise, we shall
say that an invertible $2\times 2$ matrix $M$ is an
automorphism of the binary form $G$ if $G(M\x)=G(\x)$ identically in
$\x$.  We then regard integral solutions of $G(\x)=n$ as equivalent if
and only if they are related by such an automorphism with a rational
matrix~$M$.

\specialnumber{8}\proclaim{Theorem}
Let $G(x,y)\in\Z[x,y]$ be a binary form of degree $d${\rm ,} with no factor
of multiplicity $d/2$ or more{\rm .}  Then the number of automorphisms of
$G$ is finite{\rm ,} and bounded solely in terms of $d${\rm .}  Moreover the number
of positive integers $n\leq X$ represented by the form $G$ is of exact
order $X^{2/d}${\rm ,} providing that $G(1,0)>0${\rm .}  Of these integers $n$ there are
$O_{\ep,G}(X^{52/(1+26d)+\ep})$ for which there are two or more
inequivalent integral representations{\rm .}
\endproclaim

We remark that Roth's theorem is used in the proof, so that the
implied constant is ineffective.  It seems likely, however, that this
can be avoided.

The statement that the number of representable integers is of exact
order $X^{2/d}$ is not new, and is only included for comparison with
the size of the exceptional set.  Indeed, for irreducible forms $G$, 
the lower bound is a
classical result of Erd\H{o}s and Mahler [7], dating from 1938.  In
fact \pagebreak Theorem 8 should enable one to deduce an asymptotic formula for
the number of representable integers up to $X$, such integers being
counted once only, irrespective of the number of representations.

For a form $G(x,y)=x^{e}g(x,y)$ with $e>d/2$ one can obtain
$cX^{1/(d-e)}$ representable integers merely by choosing $x=1$.  This is the
reason that such forms $G$ are excluded in the theorem.
We note also that if $G$ is a power of a quadratic form, another
excluded case, then there will be
infinitely many automorphisms, and the representations of a given integer by the form $G$ will all be equivalent.

In formulating Theorem 8 we have chosen to consider as wide a class of
forms $G$ as possible.  However for the most interesting case, in
which $G$ has no repeated factors, one can give an appreciably
stronger bound, with exponent
$$\frac{12d}{9d^{2}-6d+16}+\ep,$$
for the size of the exceptional set.  This may be achieved by using
Theorem~10 in place of Theorem 6, and taking $e=1$ in the treatment of
$S(X,C)$ in Section~7.  This remark is due to Professor Hooley.

In fact Theorem 6 does not directly entail the estimate (1.2), since the
surface $F=0$ may contain infinitely many lines.  However we may
indeed establish the following result.

\specialnumber{9}\proclaim{Theorem}
For any absolutely irreducible form
$F(\x)\in\overline{\Q}[x_{1},\ldots,x_{4}]$ of degree $d\geq 2${\rm ,}
we have 
$$N(F;B)\ll_{\ep}B^{2+\ep}.$$
\endproclaim

In higher dimensions the validity of (1.2) remains open.  We stress this
by stating formally the following conjectures.
\specialnumber{1}\proclaim{Conjecture}
For $d\geq 3$ and $n\geq 5$ we have
$$N(F;B)\ll_{\ep,F}B^{n-2+\ep}.$$
\endproclaim

\specialnumber{2}\proclaim{Conjecture}
For given $d\geq 3$ and $n\geq 5$ we have
$$N(F;B)\ll_{\ep}B^{n-2+\ep}$$
uniformly in $F$.
\endproclaim

We can do considerably better than Theorem 6 if we insist that $F$ is
nonsingular.  In this case we have the following.

\specialnumber{10}\proclaim{Theorem} \hskip-8pt
For any nonsingular form 
$F(\x)\in\overline{\Q}[x_{1},\ldots,x_{4}]$ of degree~$d${\rm ,} 
we have 
\begin{equation}
N_{1}(F;B)\ll_{\ep}B^{4/3+16/9d+\ep}.
\end{equation}
\endproclaim

For large $d$ a further improvement is possible.

\specialnumber{11}\proclaim{Theorem} \hskip-8pt
For any nonsingular form 
$F(\x)\in\overline{\Q}[x_{1},\ldots,x_{4}]$ of degree~$d${\rm ,} 
we have 
\begin{equation}
N_{1}(F;B)\ll_{\ep}B^{1+\ep}+B^{3/\sqrt{d}+2/(d-1)+\ep}.
\end{equation}
In particular
\begin{equation}
N_{1}(F;B)\ll_{\ep}B^{1+\ep},
\end{equation}
providing that $d\geq 13$.  Let $N_{2}(F;B)$ denote the number of
points counted by $N(F;B)${\rm ,} but not contained
in any curve of degree $\leq d-2$ contained in the surface{\rm .}  Then
\begin{equation}
N_{2}(F;B)\ll_{\ep}B^{3/\sqrt{d}+2/(d-1)+\ep}.
\end{equation}
Let $N_{3}(F;B)$ denote the number of
points counted by $N(F;B)${\rm ,} but not contained
in any genus zero curve of degree $\leq d-2$ contained in the surface{\rm .}  Then
\begin{equation}
N_{3}(F;B)\ll_{\ep,F}B^{3/\sqrt{d}+2/(d-1)+\ep}.
\end{equation}
\endproclaim

Thus (1.14) shows that (1.11) holds for $d\geq 13$, when $F$ is nonsingular.  

The significance of curves of degree at most $d-2$ lying in the
surface, is due to the following crucial result,
due to Colliot-Th\'{e}l\`{e}ne, and proved in the appendix.

\specialnumber{12}\proclaim{Theorem}
Let $S$ be a nonsingular surface in $\Proj^{3}${\rm ,}
of degree $d${\rm .}  Then for each degree $\delta\leq d-2$ there is a 
constant $N(\delta,d)${\rm ,} independent of $S${\rm ,} such that the surface
$S$ contains
at most $N(\delta,d)$ irreducible curves of degree $\delta${\rm .}
\endproclaim

In the case $d=3$ we have the familiar fact that a nonsingular
cubic surface has $27$ lines.  We can therefore take $N(1,3)=27$.

Since Theorem 12 shows that there are 
$O_{d}(1)$ curves of degree $\leq d-2$ in the surface, the estimate
(1.16) may be interpreted as saying that, apart from a very small number
of exceptions, all points lie on a finite number of curves of 
genus zero in the surface.

We remark that (1.13) improves on (1.12) as soon as $d\geq 6$, so that it is
only the cases $d=3,4$ and $5$ of Theorem 10 which are of real interest.
It is possible to improve the exponent $3/\sqrt{d}+2/(d-1)$ slightly,
but we shall not go into this.

There has been much work done for the special surfaces
$$F(\x)=x_{1}^{d}+x_{2}^{d}-x_{3}^{d}-x_{4}^{d}=0.$$
In particular it has been shown that for these forms $F$ we have
$$N_{1}(B)\ll B^{4/3+\ep}\;\;\;(d=3)$$
due to Heath-Brown [13],
\begin{equation}
N_{1}(B)\ll B^{5/3+\ep}\;\;\;(4\leq d\leq 7)
\end{equation}
due to Hooley [18] and [20], and
$$N_{1}(B)\ll B^{3/2+1/(d-1)+\ep}\;\;\;(d\geq 8)$$
due to Skinner and Wooley [30].  These are superseded by Theorem 11
for $d\geq 6$.  Indeed Browning, in work to appear, has shown that
(1.13) may be replaced by
$$N_{1}(F;B)\ll_{\ep}B^{2/3+\ep}+B^{3/\sqrt{d}+2/(d-1)+\ep}$$
for these particular surfaces.

For general diagonal cubic surfaces Hooley [23] showed that
$N_{1}(B)\ll_{F,{{\ep}}} B^{2}(\log B)^{-1/3}$, thereby demonstrating that
points on rational lines would dominate $N(B)$.  Moreover, also for
diagonal cubic surfaces, the author
[15] gave a conditional treatment of the bound $N_{1}(B)\ll_{F,{\ep}}
B^{3/2+\ep}$.  This is superior to Theorem 10, but assumes the Riemann
Hypothesis for the $L$-functions of elliptic curves.

We can apply our results to integral points on affine surfaces.  We
shall focus attention on the surface
$x_{1}^{d}+x_{2}^{d}+x_{3}^{d}=N$, and in view of the arithmetical
significance of this we will consider only solutions with $x_{i}>0$.
Let $r(N)$ be the number of solutions to this equation.  Then if
$d\geq 2$ we have $r(N)\ll_{d,\ep} N^{1/d+\ep}$.  No improvement in
the exponent $1/d$ has hitherto been given, for any value of $d$.  The
bound is of course best possible for $d=2$, and for $d=3$ it was shown
by Mahler [26] that $r(N)=\Omega(N^{1/12})$.  However it may be
conjectured that $r(N)\ll_{d,\ep} N^{\ep}$ as soon as $d\geq 4$.  The
mean value of $r(N)$ is also of importance.  Hua's inequality [24] shows
that
$$\sum_{n\leq B^{d}}r(n)^{2}\ll_{d,\ep} B^{7/2+\ep}$$
when $d\geq 3$.  Again, no improvement on the exponent $7/2$ has been
given hitherto, although the author [15] and Hooley [21] have shown
independently that the exponent may be reduced to $3+\ep$ in the case
$d=3$, under certain standard hypotheses concerning the Hasse-Weil
$L$-functions of cubic $3$-folds.  

We shall prove the following result.

\specialnumber{13}\proclaim{Theorem}
For $N\leq B^{d}$ we have
$$r(N)\ll_{\ep} B^{\theta+\ep}$$
where
$$\theta=\frac{2}{\sqrt{d}}+\frac{2}{d-1}.$$
It follows that
$$\sum_{n\leq B^{d}}r(n)^{2}\ll_{\ep} B^{3+\theta+\ep}.$$
\endproclaim

We note that $\theta<1$ for $d\geq 8$ and $\theta<1/2$ for $d\geq 24$.
The exponent $\theta$ may be reduced slightly with further work.

Turning to hypersurfaces of higher dimension, we have the following result.

\specialnumber{14}\proclaim{Theorem} \hskip-8pt
Let $\ep>0${\rm ,} and suppose that $B_{1},\ldots,B_{n}\geq 1$ and a form~$F${\rm ,} irreducible over $\Q${\rm ,} are given{\rm
.}  Then there exists 
$D$ depending only on $n,d$ and~$\ep${\rm ,} and an integer $k$ satisfying
$$k\ll_{\ep}
(V^{d}/T)^{d^{-(n-1)/(n-2)}}V^{\ep}(\log ||F||)^{2n-3},$$
with the following properties{\rm .}  For each $j\leq k$ there is an
integral form
$F_{j}(\x)${\rm ,} in $n$ variables{\rm ,} having degree at most $D${\rm ,} such that
\begin{itemize}
\item[{\rm 1.}] $F(\x)\nmid F_{j}(\x)$ for $1\leq j\leq k${\rm ,}
\item[{\rm 2.}] For every point $\x$ counted by $N(\b{B})$ there is an integer
$j\leq k$ such that $F_{j}(\x)=0${\rm .}
\end{itemize}

\endproclaim

Thus in particular, every point of height at most $B$ lies in 
one of at most $O_{\ep,F}(B^{(n-1)d^{-1/(n-2)}+\ep})$
proper subvarieties $F(\x)=F_{j}(\x)=0$.  The reader should note
however that such a result is trivial without a bound on the degree of
the forms $F_{j}$.  Indeed one may construct a form $F_{1}$ (with
degree dependent on $B$) such that $F_{1}(\x)=0$ for every integral
vector $\x$ in the cube $\max|x_{i}|\leq B$. 

Theorem 14 is in fact the fundamental result in this paper.  In the
case $n=3$, each point counted by $N(F;\b{B})$ lies on one of the
intersections $F(\x)=F_{j}(\x)=0$.  By B\'{e}zout's theorem, each
intersection contains at most $dD$ points, whence
$$N(F;\b{B})\ll_{\ep}(V^{d}/T)^{d^{-2}}V^{\ep}(\log ||F||)^{3}.$$
The dependence on $||F||$ can be eliminated by an appeal to Theorem 4,
so that Theorem 3 follows.

The exponents involving $1/\sqrt{d}$ appearing in our various results
all arise from the case $n=4$ of Theorem 14.  It would be remarkable if
such an exponent were optimal.  We therefore pose the following
question. 
\nonumproclaim{Question}
Is the exponent $d^{-(n-1)/(n-2)}${\rm ,} which appears in Theorem {\rm 14,} best
possible for values $n\geq 4${\rm ?}
\endproclaim

This would seem to be the single most important issue in relation to
possible sharpenings of our results.

Theorem 14 clearly opens up the prospect of results on $N(F;\b{B})$ for
$n\geq 5$.  We intend to return to this in a future paper.

This introduction would not be complete without reference to other
approaches to problems of this nature.  In particular, although the
methods developed in this paper lead in a great many cases to results
superior to those obtained hitherto, this is by no means universally
so.  The result (1.17) of Hooley is a good case in point.  Hooley uses a
sieve method, which can be thought of as counting integer vectors $\x$
for which a polynomial equation $f(X,\x)=0$ has an integral solution
$X$.  In this approach the overall number of solutions will, in 
essence, depend on the size of $\x$ alone.  In contrast, the techniques 
of the present paper produce a bound which involves the sizes both of
$X$ and~$\x$.  Thus the sieve method has potential advantages in
situations in which $X$ is large compared to $\x$. A slightly different sieve
approach, originating in work of Cohen [3], and described by the
author [10; Appendix 2], has the advantage of applying to arbitrary algebraic
hypersurfaces, but produces only $N(F;B)\ll_{\ep,F}B^{n-3/2+\ep}$.
This is inferior to the result (1.6) of Pila [27]. 
Exponential sum methods, such as those of the author [12], 
yield sharper results, but
only for nonsingular varieties.  The quality of these latter results improves
as $n$ increases.  Indeed they establish Conjecture 1, for nonsingular $F$, 
as soon as $n\geq 10$.  Other methods such as those of Schmidt [29],
depend on elementary differential geometry.  They improve slightly on
Cohen's result, and apply also to certain nonalgebraic hypersurfaces.
However none of these approaches is as effective as that of Bombieri
and Pila, for the problems considered in the present paper.

In the course of this work, the author has consulted a number of
people about issues in algebraic geometry---Jean-Louis
Colliot-Th\'{e}l\`{e}ne, Robin Hartshorne,
Miles Reid, Nick Shepherd-Barron, Sir Peter Swinnerton-Dyer, and
Yuri Tschinkel.  A number of helpful comments were also made by
Christopher Hooley.  The help of all these people is gratefully acknowledged. 

Parts of this investigation were carried out while the author was a
visitor at the Institute for Advanced Study, in Princeton.  The
hospitality and financial support of the Institute is also gratefully
acknowledged. 

\section{Preliminaries}

In this section we establish various preliminary results.

We begin by establishing Theorem 1.  This is a trivial induction
exercise.  The result is immediate for $n=2$.  In general, write 
$$F(\x)=\sum_{j=0}^{d}x_{1}^{j}F_{j}(x_{2},\ldots,x_{n}),$$
and suppose that $k$ is a value for which $F_{k}$ does not
vanish identically.  Then, by our induction assumption, there are
$O_{n,d}(B^{n-2})$ vectors $(x_{2},\ldots,x_{n})$ for which $F_{k}=0$,
and for each of these there are $O(B)$ choices for $x_{1}$.  For the
remaining vectors $(x_{2},\ldots,x_{n})$, of which there are
$O_{n}(B^{n-1})$, there are at most $d$ choices for $x_{1}$. This
produces a total of $O_{n,d}(B^{n-1})$ vectors $\x$, which completes
the induction.

We turn now to the proof of Theorem 4.  We shall write $M=\break (d+1)(d+2)/2$
and $N=d^{2}+1$, for convenience, and suppose 
that $F(\x)=0$ has solutions
$\x^{(1)},\ldots,\x^{(N)}\in Z_{n}$, where $|\x^{(i)}|\ll B$.
Consider the $N\times M$ matrix $C$, whose $i^{\rm th}$ row consists of
the $M$ possible monomials of degree $d$ in the variables 
$x^{(i)}_{1},x^{(i)}_{2},x^{(i)}_{3}$.  Then if the vector 
$\b{f}\in\Z^{M}$ has entries which are the corresponding coefficients 
of $F$, we will have $C\b{f}=\b{0}$.  Since $\b{f}\not=\b{0}$ it 
follows that $C$ has rank at most $M-1$.  Thus $C\b{g}=\b{0}$ has a
nonzero integer solution~$\b{g}$, constructed out of the
subdeterminants of $C$.  It follows that there is such a $\b{g}$ with
$|\b{g}|\ll_{d} B^{dM}$. Let $G(\x)$ be the ternary form, of degree
$d$, corresponding to the vector $\b{g}$.  Then $G(\x)$ and $F(\x)$
have at least $d^{2}+1$ common zeros, namely the vectors $\x^{(i)}$.  This will
contradict B\'{e}zout's Theorem, unless $G(\x)$ is a constant multiple of
$F(\x)$.  In the latter case $||F||\ll_{d} ||G||\ll_{d} B^{dM}$, as 
required.  This completes the proof of Theorem 4.

Many of our arguments will use elementary facts about lattices.  In
the following lemma we use $|\x|$ for the Euclidean length of the
vector $\x$.  Moreover we allow all implied constants to depend on $n$.
\specialnumber{1}\proclaim{Lemma}
\begin{itemize}
\item[{\rm (i)}] For any primitive vector $\b{c}\in\Z^{n}$ the set
$\sfLambda=\{\x\in\Z^{n}: \b{c}.\x=0\}$ is a lattice of dimension
$n-1$ and determinant $\det(\sfLambda)=|\b{c}|${\rm .}

\item[{\rm (ii)}] Let $\b{c}^{(1)},\b{c}^{(2)}\in\Z^{n}$ be nonparallel primitive
vectors{\rm ,} and let $\b{p}^{(0)}$ be the vector of length $n(n-1)/2${\rm ,} whose
coordinates are the determinants $c^{(1)}_{i}c^{(2)}_{j}-
c^{(1)}_{j}c^{(2)}_{i}${\rm ,} for $i<j$.  Write $h$ for the highest common
factor of the entries in $\b{p}^{(0)}${\rm ,} and set
$\b{p}=h^{-1}\b{p}^{(0)}${\rm .}  Then the set $\sfLambda=\{\x\in\Z^{n}:
\x\in \langle\b{c}^{(1)},\b{c}^{(2)}\rangle\}$ {\rm (}\/where $\langle\b{c}^{(1)},\b{c}^{(2)}\rangle$
denotes the $\Q$\/{\rm -}\/vector space generated by $\b{c}^{(1)}$ and $\b{c}^{(2)}${\rm )}
is a lattice of dimension $2$ and
determinant $\det(\sfLambda)=|\b{p}|${\rm .}

\item[{\rm (iii)}] Let $\sfLambda\subseteq\Z^{n}$ be a lattice of dimension $m$.
Then $\sfLambda$ has a basis $\b{b}^{(1)},\ldots,\b{b}^{(m)}$ such that if
one writes $\x\in\sfLambda$ as $\x=\sum_{j}\lambda_{j}\b{b}^{(j)}${\rm ,} then
\begin{equation}
\lambda_{j}\ll  |\x|/|\b{b}^{(j)}|.
\end{equation}
Moreover one has
\begin{equation}
\det(\sfLambda)\ll  \prod_{j=1}^{m}|\b{b}^{(j)}|\ll \det(\sfLambda).
\end{equation}

\item[{\rm (iv)}]  
Let $\x\in Z_{n}$ lie in the cube $|x_{i}|\leq B.$  Then there is a 
primitive vector $\b{y}\in\Z^{n},$ for which $\x.\b{y}=0,$ and 
such that $|\b{y}|\ll B^{1/(n-1)}${\rm .}

\item[{\rm (v)}] Let $\sfLambda\subseteq\Z^{n}$ be a lattice of dimension $m${\rm .}
Then the sphere $|\x|\leq R$ contains $O (R^{m}/\det(\sfLambda))$
points of $\sfLambda${\rm ,} providing that $R\gg\det(\sfLambda)${\rm .}

\item[{\rm (vi)}] Let $\sfLambda\subseteq\Z^{n}$ be a lattice of dimension $2${\rm .}
Then the sphere $|\x|\leq R$ contains $O(1+R^{2}/\det(\sfLambda))$
primitive points of $\sfLambda${\rm .}

\item[{\rm (vii)}] Let $P\subset\R^{2}$ be a parallelogram{\rm ,} centred on the origin{\rm ,}
having area $A${\rm .}  Then $P$ contains $O(1+A)$ primitive integer vectors{\rm .}
\end{itemize}

\endproclaim

Statement (i) of the lemma is a special case of Heath-Brown [11; Lemma~1].

For part (ii), we first note that it is trivial that $\sfLambda$ is a
two-dimensional lattice.  Choose a basis $\b{b}^{(1)},\b{b}^{(2)}$ for
$\sfLambda$, and set
$\b{c}^{(i)}=a_{i1}\b{b}^{(1)}+a_{i2}\b{b}^{(2)}$, for $i=1,2$.  If
$\b{q}^{(0)}$ is the vector formed from the determinants 
$b^{(1)}_{i}b^{(2)}_{j}-b^{(1)}_{j}b^{(2)}_{i}$ for $i<j$, 
then $|\b{q}^{(0)}|$
is the area of the parallelogram spanned by $\b{b}^{(1)}$ and
$\b{b}^{(2)}$, so that $|\b{q}^{(0)}|=\det(\sfLambda)$.  Moreover, we
will have $\b{p}^{(0)}=(a_{11}a_{22}-a_{12}a_{21})\b{q}^{(0)}$, so
that $\b{p}$ will be a scalar multiple of $\b{q}^{(0)}$.  To complete
the proof of part (ii) it therefore suffices to show that
$\b{q}^{(0)}$ is primitive.  However if $p$ were a prime dividing
$\b{q}^{(0)}$ then the reductions modulo $p$ of $\b{b}^{(1)}$ and
$\b{b}^{(2)}$ would be proportional.  There would then be integers
$\lambda_{1},\lambda_{2}$, not both multiples of $p$, and an
integral vector $\b{b}$, such that $\lambda_{1}\b{b}^{(1)}+
\lambda_{2}\b{b}^{(2)}=p\b{b}$.  Since we then have $\b{b}\in\sfLambda$,
this would contradict the fact that $\b{b}^{(1)},\b{b}^{(2)}$ is a
basis for $\sfLambda$.

For statement (iii)  we note that
Davenport [4; Lemma 5]  shows the existence of
a basis $\b{b}^{(j)}$ with the property (2.1).  Moreover in the course of
the proof he shows [4; (14)] that
$$\prod_{j=1}^{m}|\b{b}^{(j)}|\ll \det(\sfLambda).$$
It is of course trivial that
$$\prod_{j=1}^{m}|\b{b}^{(j)}|\gg \det(\sfLambda),$$
for any basis.  

For part (iv) we note that the lattice of integral
vectors $\b{y}$ satisfying\break $\x.\b{y}=0$ has dimension $n-1$ and 
determinant $|\b{x}|$, by
statement (i).  According to part (iii) there is therefore a basis 
element $\b{y}'$, say, with $|\b{y}'|\ll |\x|^{1/(n-1)}$, which is
sufficient. 

Since the condition $R\gg\det(\sfLambda)$ ensures that the basis
vectors in part (iii) all satisfy $|\b{b}^{(j)}|\ll R$, the
fundamental parallelpiped formed from these will fit inside a suitable
constant multiple of the sphere $|\x|\leq R$.  
Statement (v) of the lemma then follows.

To establish part (vi) we note that $\sfLambda$ has a basis
$\b{b}^{(1)},\b{b}^{(2)}$ as in part (iii).  Thus if $\x=
\lambda_{1}\b{b}^{(1)}+\lambda_{2}\b{b}^{(2)}$ satisfies $|\x|\leq R$,
then $\lambda_{j}\ll R/|\b{b}^{(j)}|$, for $j=1,2$.  There are
therefore 
$$\ll R^{2}|\b{b}^{(1)}|^{-1}|\b{b}^{(2)}|^{-1}\ll
R^{2}\det(\sfLambda)^{-1}$$
possible pairs $\lambda_{1},\lambda_{2}$ with
$\lambda_{1}\lambda_{2}\not=0$.  Moreover, since $\x$ is to be
primitive, we can only have $\x=\pm\b{b}^{(1)}$ or $\pm\b{b}^{(2)}$
when $\lambda_{1}\lambda_{2}=0$.  This suffices for part (vi).

For the final assertion, we begin by constructing a rectangle $P'$ 
including~$P$, 
centred on the origin, and having area $A'\ll A$.  We may then produce
an ellipse $E$ centred on the origin, and having area $A''\ll A'\ll
A$.  The desired estimate is then a corollary of Heath-Brown [11;
Lemma 2].

We shall also want some results from elimination theory.  We first
state without proof the following basic result.

\specialnumber{2}\proclaim{Lemma}
Let integers $m\geq n\geq 2$ and $d\geq 1$ be given{\rm .}  Then there exist
integers $m',d'$ depending at most on $m$ and $d${\rm ,} as follows{\rm .}  Let
$F_{1}(\x),\ldots,F_{m}(\x)$ be forms in $n$ variables{\rm ,} with
coefficients in $\Qb${\rm ,} each with
degree at most $d${\rm .}  Let $\b{f}_{i}$ be the coefficient vector of
$F_{i}${\rm .} Then there exist polynomials
$E_{i}(\b{f}_{1},\ldots,\b{f}_{m})$ for $1\leq i\leq m'$
over $\Qb${\rm ,} with 
the following properties{\rm .}
\begin{itemize}
\item[{\rm 1.}] Each $E_{i}$ has total degree at most $d'${\rm .}
\item[{\rm 2.}] The polynomials $E_{i}$ are
homogeneous functions in each $\b{f}_{j}${\rm .} 
\item[{\rm 3.}] The simultaneous equations
$F_{i}(\x)=0${\rm ,} for $1\leq i\leq m${\rm ,}
have a nonzero solution over $\Qb$ if and only if
$E_{i}(\b{f}_{1},\ldots,\b{f}_{m})=0$ for $1\leq i\leq m'${\rm .}
\end{itemize}
\endproclaim

Note that the lemma does not assert that the $E_{i}$ are nonzero.

From this we shall deduce the following.

\specialnumber{3}\proclaim{Lemma}
Let $F(\x)\in\Qb[x_{1},x_{2},x_{3},x_{4}]$ be a form of degree $d${\rm .}
For any $\b{y}\not=\b{0}${\rm ,} let $H_{{\bf y}}(x_{i},x_{j},x_{k})$ be
a form got by eliminating 
one variable from the equations $F(\x)=\b{y}.\x=0${\rm .}
Then for any positive integer $\delta<d$ there is an integer
$m=O_{d}(1)${\rm ,} and forms 
$E_{i,\delta}(\b{y})\in\Qb[y_{1},y_{2},y_{3},y_{4}]${\rm ,} with $1\leq i
\leq m${\rm ,} whose degrees are
bounded in terms of $d${\rm ,} and which vanish simultaneously
precisely at those points $\b{y}\not=0$ for which $H_{{\bf y}}$ has a
factor of degree $\delta${\rm .}
\endproclaim

Again we do not assert that the forms $E_{i,\delta}$ are nonzero.  We
note that it does not matter how we eliminate one of the variables to
produce $H_{\bf y}$, since if one of the resulting forms has a factor of
degree $d$ they all will.

To deduce the above result from Lemma 2 we consider possible factors
of degree $\delta$ in 
$$y_{1}^{d}F(\x)=F(-y_{2}x_{2}-y_{3}x_{3}-y_{4}x_{4},y_{1}x_{2},
y_{1}x_{3},y_{1}x_{4})= f(x_{2},x_{3},x_{4}),$$
say.  If $G$ is a nonzero form of degree $\delta$, the relation 
\begin{equation}
\lambda f(x_{2},x_{3},x_{4})=G(x_{2},x_{3},x_{4})H(x_{2},x_{3},x_{4})
\end{equation}
produces a system $L_{i}(\lambda,\b{h})=0$ of homogeneous linear
equations in $\lambda$ and the coefficients $\b{h}$, say, of $H$.  
The coefficients
of the $L_{i}$ are polynomials in the $y_{i}$ and in the coefficients
$\b{g}$, say,
of $G$.  According to Lemma 2 we produce polynomials $E_{j}(\b{y},\b{g})$
in these latter variables, which vanish precisely when (2.3) has a
nonzero solution.  (In this case Lemma 2 is a well-known result in
linear algebra, the polynomials $E_{j}$ arising as determinants.)
If $\lambda$ were to vanish in such a solution,
then the form $H$ must vanish too, since a polynomial ring over a
field has no zero-divisors.  Thus, if $G$ is nonzero, then $G$ divides
$f$ precisely when the polynomials $E_{j}(\b{y},\b{g})$ all vanish.  Since
the divisibility of $f$ by $G$ is unnaffected by replacing 
$G$ by $cG$, or $\b{y}$ by
$c'\b{y}$, for any nonzero $c,c'$, we see that the various
bi-homogeneous parts $B_{k}(\b{y},\b{g})$ of $E_{j}(\b{y},\b{g})$ must vanish 
precisely when $G|f$.  We note that the degrees of the forms $B_{k}$,
and the number of forms that arise, are $O_{d}(1)$.

A second application of Lemma 2 now produces forms $J_{l}(\b{y})$,
which vanish simultaneously if and only if there is a nonzero set of
coefficients $\b{g}$ for $G$ which make all the forms $B_{k}(\b{y},\b{g})$
vanish.  As we have seen, this is precisely equivalent to the
requirement that $f$ should have a factor of degree $\delta$.  We
rename the forms $J_{l}(\b{y})$ as $J_{l,1}(\b{y})$, to denote the
fact that $x_{1}$ was eliminated in forming $f$.  Thus, in a
precisely analogous way, we produce forms $J_{l,i}(\b{y})$ for
$i=2,3,4$.  If $H_{\b{y}}$ has a factor of degree $\delta$, then all
four of the possible forms $f$ have such a factor, so that 
$J_{l,i}(\b{y})=0$ for each $l$ and each $i$.  Conversely this latter
condition implies that each of the four possible forms $f$ factors.
Since at least one of the $y_{i}$ is nonzero, this implies that
$H_{\bf y}$ has a factor of degree $d$.  We may therefore take the
forms $E_{i,\delta}$ to be the various $J_{l,i}$.  Clearly both 
the number of such forms, and their degrees, are bounded in
terms of $d$.

\section{Proof of Theorem 14}

Before beginning the proof of the above theorem we shall
require a preliminary result.
Let
$$S(F;\b{B},p)=\{\x\in Z_{n}:\,F(\x)=0,\,
|x_{i}|\leq B_{i},\,(1\leq i\leq 1),\,p\nmid\nabla F(\x)\},$$
and
$$S(F;\b{B})=\{\x\in Z_{n}:\,F(\x)=0,\,
|x_{i}|\leq B_{i},\,(1\leq i\leq 1),\,\nabla F(\x)\not=\b{0}\}.$$
We then have the following lemma.

\specialnumber{4}\proclaim{Lemma}
Let $B=2+\max|B_{i}|$ and $r=[\log(||F||B)]${\rm ,} and suppose that 
$$P\geq \log^{2}(||F||B).$$
Then there are distinct primes
$p_{1},\ldots,p_{r}${\rm ,} such that $P\ll p_{i}\ll P$ and 
$$S(F;\b{B})=\bigcup_{i=1}^{r}S(F;\b{B},p_{i}).$$
\endproclaim

We remark that this result is the sole point
at which a dependence on $||F||$ enters our arguments.  To prove Lemma
4 we merely choose the primes $p_{i}$ as the first $r$ primes
$p_{i}>AP$, for a suitable constant $A$.  
Since $P\gg r^{2}$ this yields $P\ll p_{i}\ll P$.  Now
if $\x$ is in $S(F;\b{B})$, then some partial derivative $\partial
F/\partial x_{j}$, say, must be nonzero.   Since 
$$\frac{\partial F}{\partial x_{j}}\ll_{n} ||F||B^{d-1},$$
it follows that
$$\#\left\{p>AP:\,p|\frac{\partial F}{\partial x_{j}}\right\}\ll_{n,d} 
\frac{\log||F||B}{\log AP}.$$
Thus there are fewer than $r$ such primes, if $A$ is large enough.  We
therefore see that there is some prime $p_{i}$ which does not divide
$\partial F/\partial x_{j}$, whence $\x\in S(F;\b{B},p_{i})$, as required 
for the lemma.

To prove Theorem 14 we begin by considering singular points.
Any singular points of $F(\x)=0$ satisfy
$$\frac{\partial F(\x)}{\partial x_{i}}=0,\;\;\; (1\leq i\leq n).$$
Since $F$ is irreducible, at least one of the forms $\partial
F/\partial x_{i}$ is not identically zero.  Clearly such a form cannot be a
multiple of $F$ since its degree is $d-1$.  We therefore include 
one of the partial derivatives of $F$ amongst the forms $F_{i}$
described in Theorem 14 to take care of the singular points of $F(\x)$.

It therefore remains to examine nonsingular points, and here we apply
Lemma 4.  This shows that we may consider points that are nonsingular
modulo a suitable prime $p$, at a cost of a factor $\log(||F||B)$ in
our final estimate for $k$.

With this understanding, we now define $k$ to be the number of
nonsingular points $\b{t}\in\Proj^{n-1}(\F_{p})$ on $F(\b{t})=0$.  Thus
$k\ll p^{n-2}\ll P^{n-2}$, and we split the points $\x\in 
S(F;\b{B},p)$ into $k$ sets 
$$S(\b{t})=\{\x\in S(F;\b{B},p):
\x\equiv\rho\b{t}\mod{p}\;\mbox{for some}\;\rho\in\Z\}.$$
Our aim is
to show that if $P$ is chosen so that
\begin{equation}
P\gg (V^{d}/T)^{(n-2)^{-1}d^{-(n-1)/(n-2)}}V^{\ep}\log^{2}||F||,
\end{equation}
then for each set $S(\b{t})$
there is a corresponding form $F_{j}$ such that $F_{j}(\x)=0$ for all
$\x\in S(\b{t})$.  Note that the term $\log^{2}||F||$ has been
included above so as to ensure that $P$ is acceptable for Lemma 4.

From now on we shall focus our attention on a fixed $\b{t}$.  
Since $\b{t}$ is nonzero we
may suppose without loss of generality that $t_{i}=1$ for some $i$,
and, again without loss of generality, we may take $i=1$.  If
\begin{equation}
\frac{\partial F}{\partial x_{i}}(\b{t})=0,\;\;\;(2\leq i\leq n)
\end{equation}
then 
$$0=dF(\b{t})=\b{t}.\nabla F(\b{t})=
\frac{\partial F}{\partial x_{1}}(\b{t}),$$
whence $\nabla F(\b{t})=\b{0}$.  This contradiction shows that one of
the partial derivatives in (3.2) must be nonvanishing, and we assume,
without loss of generality, that
\begin{equation}
\frac{\partial F}{\partial x_{2}}(\b{t})\not=0.
\end{equation}

We proceed
to lift $\b{t}$ to a $p$-adic solution $\b{u}\in\Z_{p}^{n}$ of
$F(\b{u})=0$.  In view of (3.3), Hensel's lemma may be used to produce a
solution in which $u_{1}=1$.  We now require the following result.

\specialnumber{5}\proclaim{Lemma}
Let $F(\x)\in\Z_{p}[\x]$ be a form in $n$ variables{\rm ,} and suppose that
$\b{u}\in\Z_{p}^{n}$ satisfies $u_{1}=1$ and
$$F(\b{u})=0,\;\;\; p\nmid\frac{\partial F}{\partial x_{2}}(\b{u}).$$
Then{\rm ,} for any integer $m\geq 1$ there exists
$f_{m}(Y_{3},Y_{4},\ldots,Y_{n})\in\Z_{p}[Y_{3},\ldots,Y_{n}]${\rm ,} such
that if $F(\b{v})=0$ for some $\b{v}\in\Z_{p}^{n}$ with $v_{1}=1$ and
$\b{v}\equiv\b{u}\pmod{p}${\rm ,} then
\begin{equation}
v_{2}\equiv f_{m}(v_{3},\ldots,v_{n})\mod{p^{m}}.
\end{equation}
\endproclaim

One could alternatively formulate Lemma 5 to say that, for given
$v_{3},\ldots,v_{n}$, the equation $F(\b{v})=0$ has a unique solution
$v_{2}$, and that this solution is given by a $\Z_{p}$-integral
power series $v_{2}=
f(v_{3},\ldots,v_{n})$.  One could then use such a result in what
follows to replace the sequence of polynomials $f_{m}$.

For the proof of Lemma 5 let
$$\frac{\partial F}{\partial x_{2}}(\b{u})=\mu,$$
say, and define the polynomials $f_{m}$ inductively by taking
$f_{1}(Y_{3},\ldots,Y_{n})=u_{2}$, (constant) and
$$f_{m+1}(Y_{3},\ldots,Y_{n})=f_{m}(Y_{3},\ldots,Y_{n})-\mu^{-1}
F(1,f_{m}(Y_{3},\ldots,Y_{n}),Y_{3},\ldots,Y_{n}),$$
for $m\geq 1$.  Clearly Lemma 5 now holds for $m=1$.  We prove the general
case by induction on $m$.  Thus we may suppose that
$$v_{2}\equiv f_{m}(v_{3},\ldots,v_{n})\mod{p^{m}},$$
and we write
$$v_{2}=f_{m}(v_{3},\ldots,v_{n})+\lambda p^{m},$$
where $\lambda\in\Z_{p}$.  Then
\begin{eqnarray}
0&=&F(\b{v})\\
&\equiv &F\lrp{1,f_{m}(v_{3},\ldots,v_{n}),v_{3},\ldots,v_{n}}\nonumber\\
&& +\ \lambda p^{m}\frac{\partial F}{\partial x_{2}}
\lrp{1,f_{m}(v_{3},\ldots,v_{n}),v_{3},\ldots,v_{n}}\mod{p^{m+1}}.\nonumber
\end{eqnarray}
Since $\b{v}\equiv\b{u}\pmod{p}$, the induction hypothesis (3.4) shows that
$$f_{m}(v_{3},\ldots,v_{n})\equiv u_{2}\mod{p},$$
and hence that
$$\frac{\partial F}{\partial x_{2}}
(1,f_{m}(v_{3},\ldots,v_{n}),v_{3},\ldots,v_{n})\equiv\mu\mod{p}.$$
The congruence (3.5) then implies
$$\lambda p^{m}\equiv
-\mu^{-1}F(1,f_{m}(v_{3},\ldots,v_{n}),v_{3},\ldots,v_{n})\mod{p^{m+1}},$$
whence
$$v_{2}\equiv f_{m+1}(v_{3},\ldots,v_{n})\mod{p^{m+1}},$$
as required for the induction step.

We are now ready to examine our set $S(\b{t})$.  Let $\x\in S(\b{t})$,
so that the reduction modulo $p$ of $\x$ represents the same
projective point as does $\b{t}$.  Thus $p\nmid x_{1}$ so that we may
interpret $x_{1}^{-1}\x=\b{v}$, say, as a vector in $\Z_{p}^{n}$.  We then
see that $v_{1}=1$ and $v_{i}=u_{i}+y_{i}$ for $2\leq i\leq n$, for
suitable $y_{i}\in p\Z_{p}$. We shall define a collection of
monomials of degree $D$, by choosing a set of exponents 
$${\cal E}\subseteq\left\{(e_{1},\ldots,e_{n})\in\Z^{n}:\,e_{i}\geq
0,\,(1\leq i\leq n),\;\sum_{i=1}^{n}e_{i}=D\right\},$$
and considering monomials of the form
$$X_{1}^{e_{1}}\ldots X_{n}^{e_{n}}=\b{X}^{\bf e},$$
say.  We shall write $E=\#{\cal E}$, and suppose that $E\leq \#S(\b{t})$.
Now take distinct elements $\x^{(1)},\ldots,\x^{(E)}$ of $S(\b{t})$ and
consider the $E\times E$ determinant
$$\Delta=\det(\x^{(i)\b{e}})_{1\leq i\leq E,\;\b{e}\in\cl{E}},$$
with rows corresponding to the different vectors $\x^{(i)}$ and columns
corresponding to the various exponent $n$-tuples $\b{e}$.  Our first task is
to show that $\Delta$ must vanish, if $p$ is sufficiently large in
terms of the various $B_{i}$.

We begin by considering $\Delta$ modulo a large power $p^{m}$ of $p$.
We have
$$\Delta=\Big(\prod_{1\leq i\leq E}x_{1}^{(i)}\Big)^{D}
\det(\b{v}^{(i)\b{e}})_{1\leq i\leq E,\;\b{e}\in\cl{E}},$$
with $\b{v}^{(i)}=(x_{1}^{(i)})^{-1}\x^{(i)}$, as above.  According to
Lemma 5 we deduce that
$$\Delta\equiv\Big(\prod_{1\leq i\leq E}x_{1}^{(i)}\Big)^{D}\Delta_{0}
\mod{p^{m}},$$
where
$$\Delta_{0}=\det(M_{0}),\;\;\;
M_{0}=(\b{w}^{(i)\b{e}})_{1\leq i\leq E,\;\b{e}\in\cl{E}},$$
with
$$w_{1}^{(i)}=1,\;\;\;w_{2}^{(i)}=f_{m}\lrp{v_{3}^{(i)},\ldots,v_{n}^{(i)}},$$
and
$$w_{j}^{(i)}=v_{j}^{(i)}\;\;\;(3\leq j\leq n).$$
We now set $v_{j}^{(i)}=u_{j}+y_{j}^{(i)}$ for $3\leq j\leq n$, so
that $p|y_{j}^{(i)}$.  Thus
$$\b{w}^{(i)\b{e}}=w_{1}^{(i)e_{1}}\ldots w_{n}^{(i)e_{n}}=
g_{\bf e}\lrp{y_{3}^{(i)},y_{4}^{(i)},\ldots,y_{n}^{(i)}}$$
for an appropriate set of polynomials $g_{\bf e}(Y_{3},\ldots,Y_{n})
\in\Z_{p}[Y_{3},\ldots,Y_{n}]$.  We now introduce an ordering on the
exponent vectors 
$$\b{f}=(f_{3},\ldots,f_{n}),\;\;\;(f_{j}\in\Z,\;\;f_{j}\geq 0),$$
by setting $\b{f}\prec\b{f}'$ if either 
\begin{itemize}
\item[1.] $\sum f_{j}<\sum f_{j}'$, or 
\item[2.] $\sum f_{j}=\sum f_{j}'$, and there is some $j$ such that
$f_{h}=f_{h}'$ for $h<j$ but $f_{j}<f_{j}'$.
\end{itemize}
As the reader will observe, it is important, in what follows, to have
$\b{f}\prec\b{f}'$ in case 1, but the ordering when $\sum f_{j}=\sum
f_{j}'$ is immaterial.  We shall order the monomials $\b{Y}^{\b{f}}$
in the analogous way.

We proceed to perform column operations on $M_{0}$ as follows.
We look for the `smallest' monomial $\b{Y}^{\bf f}$, say, occurring 
in any of the polynomials 
$g_{\bf e}$.  If this monomial occurs in more than
one such polynomial we take the occurrence for which the coefficient has
the smallest $p$-adic order.  We swap columns to bring this term into
the first column, and then subtract $p$-adic integer multiples of the
new first column from all those columns containing the monomial 
$\b{Y}^{\bf f}$, so as to remove it entirely, except from the first
column.  This process is then repeated with the remaining $n-1$
columns, looking again for the `smallest' monomial, moving it to
column $2$ and removing it from all subsequent columns.  We proceed in
this way to obtain an expression
$$\Delta_{0}=\det(M_{1}),\;\;\;M_{1}=
\big(h_{e}(y_{3}^{(i)},\ldots,y_{n}^{(i)})\big)
_{1\leq i\leq E,\;1\leq e\leq E},$$
in which one has polynomials $h_{e}(\b{Y})\in\Z_{p}[\b{Y}]$, with
successively larger `smallest' monomial terms.  The number of monomials of
total degree $f$ is 
$$\left(\begin{array}{cc} f+n-3\\ n-3\end{array}\right)=n(f),$$
say.  Thus if $e>n(0)+n(1)+\ldots+n(f-1)$, the `smallest' term
in $h_{e}(\b{Y})$ must have total degree at least $f$.  Since
$p|y_{j}^{(i)}$ for $3\leq j\leq n$ we deduce that every element in the
$e^{\rm th}$ column of $M_{1}$ must be divisible by $p^{f}$.  We note that
$$\sum_{i=0}^{f}n(i)=\left(\begin{array}{cc} f+n-2\\ n-2\end{array}\right),$$
and that
$$\sum_{i=0}^{f}in(i)=(f+1)\left(\begin{array}{cc} 
f+n-2\\ n-2\end{array}\right)-
\left(\begin{array}{cc} f+n-1\\ n-1\end{array}\right).$$
It therefore follows that if 
\begin{equation}
\left(\begin{array}{cc} f+n-2\\ n-2\end{array}\right)\leq E<
\left(\begin{array}{cc} (f+1)+n-2\\ n-2\end{array}\right),
\end{equation}
then $\Delta_{0}$ is divisible by 
$$p^{n(1)+2n(2)+\ldots+fn(f)+(f+1)(E-n(0)-n(1)-\ldots-n(f))}=p^{\nu},$$
say, where
\begin{equation}
\nu=(f+1)E-\left(\begin{array}{cc} f+n-1\\ n-1\end{array}\right).
\end{equation}
If we choose our original prime power $p^{m}$ to have $m=\nu$ we may
therefore conclude as follows.

\specialnumber{6}\proclaim{Lemma}
Let $E$ lie in the range {\rm (3.6),} and suppose that $\nu$ is given by {\rm (3.7).}
Then
$$\nu_{p}(\Delta)\geq\nu.$$
\endproclaim

We shall compare this result with information on the size of
$\Delta$.  Since $|x_{j}^{(i)}|\leq B_{j}$, every element of the
column corresponding to exponent vector $\b{e}$ has modulus at most
$\b{B}^{\bf e}$. Thus an elementary estimate yields
$$|\Delta|\leq E^{E}\prod_{{\bf e}\in\cl{E}}\b{B}^{\bf e}.$$
We shall set
\begin{equation}
\sum_{{\bf e}\in\cl{E}}\b{e}=\b{E},
\end{equation}
say, and require that
\begin{equation}
p^{\nu}>E^{E}\b{B}^{\bf E}.
\end{equation}
Then $\Delta$ must vanish.

In forming $\Delta$ we assumed that $\#S(\b{t})\geq E$, and we
took $\x^{(1)},\ldots,\x^{(E)}$ to be any distinct elements of 
$S(\b{t})$.  Thus if we set $\#S(\b{t})=K$ and
consider the matrix
$$M_{2}=\big(\x^{(i)\b{e}}\big)_{1\leq i\leq K,\;\b{e}\in\cl{E}},$$
where $\x^{(i)}$ now runs over all elements of $S(\b{t})$, we see
that $M_{2}$ can have rank at most $E-1$.  This is trivial when $K\leq
E-1$, and otherwise every $E\times E$ minor vanishes, by what we have
proved. It follows that $M_{2}\b{c}=\b{0}$ for some nonzero vector
$\b{c}\in\Z^{E}$. Thus, if we set
\begin{equation}
G(\b{X})=\sum_{{\bf e}\in\cl{E}}c_{\bf e}\b{X}^{\bf e},
\end{equation}
we have produced a nonzero polynomial, of degree $D$, and such 
that $G(\x)=0$ for every $\x\in S(\b{t})$.

It remains to select the exponent set $\cl{E}$ so as to ensure that
$F(\x)\nmid G(\x)$.  We write
$$F(X_{1},\ldots,X_{n})=
\sum_{{\bf f}}a_{{\bf f}}X_{1}^{f_{1}}\ldots X_{n}^{f_{n}},$$
and consider the Newton polyhedron $P$, defined as the convex hull of the
points $\b{f}\in\R^{n}$ for which $a_{{\bf f}}\not=0$.  Clearly $P$ is
a subset of the affine hyperplane given by $\sum f_{i}=d$.  Any vertex 
of $P$ will be an exponent vector $\b{f}$, with $a_{{\bf f}}\not=0$.
Now consider such a vertex $\b{f}^{*}$, say, at which
$$\sum_{i=1}^{n} f^{*}_{i}\log B_{i}$$
is maximal.  We proceed to choose numbers $B_{i}'$ in the range
$[B_{i}\,,\,1+B_{i}]$, such that the values of $\log B_{i}'$ are linearly
independent over $\Q$, and such that
$$\sum_{i=1}^{n} f_{i}\log B_{i}'$$
is maximal only at the vertex $\b{f}^{*}$ of $P$.  Let
the maximal value be $M_{F}$.

Suppose now that $G(\b{X})$ is given by (3.10), and that $G(\x)$ is a
multiple of $F(\b{X})$, so that $G(\b{X})=F(\b{X})K(\b{X})$, say.  
Let
$$K(X_{1},\ldots,X_{n})=
\sum_{{\bf k}\in\cl{K}}b_{{\bf k}}X_{1}^{k_{1}}\ldots X_{n}^{k_{n}},$$
with $b_{{\bf k}}\not=0$, and suppose that
$$\sum_{i=1}^{n} k_{i}\log B_{i}'$$
is maximal at $\b{k}=\b{k}^{*}$, say, with maximal value $M_{K}$.  
Clearly $\b{k}^{*}$ is
unique, since the $\log B_{i}'$ are linearly independent over $\Q$.
Now all terms
$$a_{{\bf f}}X_{1}^{f_{1}}\ldots X_{n}^{f_{n}}.
b_{{\bf k}}X_{1}^{k_{1}}\ldots X_{n}^{k_{n}}$$
arising from the product $F(\b{X})K(\b{X})$ will have
$$\sum_{i=1}^{n} (f_{i}+k_{i})\log B_{i}'<M_{F}+M_{K},$$\pagebreak

\noindent
with the exception of the term for $\b{f}=\b{f}^{*}$ and
$\b{k}=\b{k}^{*}$.  It follows that the monomial
$$X_{1}^{f^{*}_{1}+k^{*}_{1}}\ldots X_{n}^{f^{*}_{n}+k^{*}_{n}}$$
occurs in $G(\x)$ with nonzero coefficient.  

We now define 
$$\cl{E}=\left\{(e_{1},\ldots,e_{n})\in\Z^{n}:\,e_{i}\geq
0,\,(1\leq i\leq
n),\;\sum_{i=1}^{n}e_{i}=D,\;e_{i}<f^{*}_{i}\;\mbox{for some}\;i\right\}.$$
In the light of the above discussion it is then apparent that 
we cannot have $F(\b{X})|G(\b{X})$.

It remains to choose the parameter $D$.  We see from (3.9) that it
suffices to require that
$$p\gg_{D}\;\prod_{i=1}^{n} B_{i}^{E_{i}/\nu}.$$
However it is an elementary matter to calculate that if $D\geq d$ then
$$E=\left(\begin{array}{cc} D+n-1\\ n-1\end{array}\right)-
\left(\begin{array}{cc} D-d+n-1\\ n-1\end{array}\right)
=\frac{dD^{n-2}}{(n-2)!}+O(D^{n-3}).$$
Here we follow the convention that implied constants may depend on $n$
and $d$.
Moreover, since (3.6) implies that 
$$E=\frac{f^{n-2}}{(n-2)!}+O(f^{n-3}),$$
we deduce that
$$f=d^{1/(n-2)}D+O(1).$$
Thus (3.7) yields
\begin{eqnarray}
\nu &=& \frac{(n-2)f^{n-1}}{(n-1)!}+O(f^{n-2})\\
&=& d^{(n-1)/(n-2)}(n-2)\frac{D^{n-1}}{(n-1)!}+O(D^{n-2}).\nonumber
\end{eqnarray}
In order to find the vector $\b{E}$ defined in (3.8), we 
write $\cl{E}=\cl{E}_{1}\setminus\cl{E}_{2}$, where
$$\cl{E}_{1}=\left\{(e_{1},\ldots,e_{n})\in\Z^{n}:\,e_{i}\geq
0,\,(1\leq i\leq
n),\;\sum_{i=1}^{n}e_{i}=D\right\}$$
and
$$\cl{E}_{2}=\left\{(e_{1},\ldots,e_{n})\in\Z^{n}:\,e_{i}\geq
0,\,(1\leq i\leq
n),\;\sum_{i=1}^{n}e_{i}=D,\;e_{i}\geq f^{*}_{i}\;\mbox{for all}\;i\right\}.$$
Then
\begin{eqnarray*}
\sum_{{\bf e}\in\cl{E}_{1}}e_{i}&=&
\frac{1}{n}\sum_{{\bf e}\in\cl{E}_{1}}\sum_{i=1}^{n}e_{i}\\[4pt]
&=&\frac{D}{n}\#\cl{E}_{1}\\[4pt]
&=&\frac{D}{n}\left(\begin{array}{cc} D+n-1\\[4pt] n-1\end{array}\right),
\end{eqnarray*}
and similarly,
\begin{eqnarray*}
\sum_{{\bf e}\in\cl{E}_{2}}e_{i}&=&\lrp{f^{*}_{i}+\frac{D-d}{n}}\#\cl{E}_{2}\\[4pt]
&=&\lrp{f^{*}_{i}+\frac{D-d}{n}}
\left(\begin{array}{cc} D-d+n-1\\[4pt] n-1\end{array}\right).
\end{eqnarray*}
Thus
\begin{eqnarray*}
E_{i}&=&\sum_{{\bf e}\in\cl{E}}e_{i}\\[4pt]
&=&\frac{D}{n}\left(\begin{array}{cc} D+n-1\\[4pt] n-1\end{array}\right)-
\lrp{f^{*}_{i}+\frac{D-d}{n}}
\left(\begin{array}{cc} D-d+n-1\\[4pt] n-1\end{array}\right)\\[4pt]
&=&(d-f^{*}_{i})\frac{D^{n-1}}{(n-1)!}+O(D^{n-2}).
\end{eqnarray*}
In view of (3.11) we find that
$$E_{i}/\nu=(n-2)^{-1}(d-f^{*}_{i})d^{-(n-1)/(n-2)}+O(D^{-1}),$$
whence it suffices to have
$$p\gg (V^{d}/T)^{(n-2)^{-1}d^{-(n-1)/(n-2)}}V^{O(1/D)}.$$
The condition (3.1) is therefore sufficient, providing that we take 
$D\geq D(n,d,\ep)$.  This completes the proof of Theorem 14.

\section{Curves in $\Proj^{3}$}

In this section we shall prove Theorem 5, by projecting the curve $C$ onto
a suitable planar curve.  
The following 
result shows how this may be done without changing the degree of the
curve.  Recall that the degree of a curve in $\Proj^{3}$
may be defined as the number
of points of intersection with a generic plane. \pagebreak

\specialnumber{7}\proclaim{Lemma}
Let $C\subset\Proj^{3}$ 
be an irreducible projective curve of degree $d${\rm .}  Then there are
nonzero integer vectors $\b{y}$ and $\b{c}$ with $|\b{y}|,|\b{c}|
\ll 1${\rm ,} such that $\b{y}.\b{c}\not=0${\rm ,} and so that the 
projection of $C$ parallel to $\b{y}${\rm ,} onto
the plane $\x.\b{c}=0$ produces an irreducible curve of degree $d${\rm .}
Moreover each fibre contains at most $d$ points{\rm .}
\endproclaim

It is a familiar fact that the generic projection of $C$ onto a plane
will indeed be an irreducible curve of degree $d$.  Thus the thrust of
the result is that we can choose a projection with $|\b{y}|\ll 1$.
One difficulty in the proof is that we do not have a convenient basis
for the ideal of polynomials vanishing on $C$.

Before proving Lemma 7, we show how Theorem 5 follows.  Write $\pi$
for the projection given by Lemma 7.  
If $\x\in Z_{4}$, with $|\x|\ll B$, then 
$$\pi(\x)=\x-\frac{(\x.\b{c})}{(\b{y}.\b{c})}\b{y},$$
whence $(\b{y}.\b{c})\pi(\x)$ is an integral vector, 
with $|(\b{y}.\b{c})\pi(\x)|\ll B$.  
Although the vectors 
$(\b{y}.\b{c})\pi(\x)$ may not be primitive, there are, according to
Lemma 7, at most $d$ values of $\x$ for which $(\b{y}.\b{c})\pi(\x)$
is projectively equivalent to a given point in the plane
$\b{z}.\b{c}=0$.  Thus it will suffice to show that
the curve $\pi(C)$ has $O_{\ep}(B^{2/d+\ep})$ points in
the region $|\b{z}|\ll  B$.

According to parts (i) and (iii) of
Lemma 1, we can choose a basis for the lattice of integer
vectors in the plane $\b{z}.\b{c}=0$, with respect to which
$\b{z}$ will have coordinates
$(\lambda_{1},\lambda_{2},\lambda_{3})$ with $\lambda_{i}\ll B$.
Since Lemma 7 produces a curve $\pi(C)$ of degree $d$, we may
apply Theorem 3 to show that the number of primitive points
$(\lambda_{1},\lambda_{2},\lambda_{3})$ on the curve $\pi(C)$,
lying in the region $\lambda_{i}\ll B$, is indeed
$O_{\ep}(B^{2/d+\ep})$.  This establishes Theorem 5.

The remainder of this section is devoted to the proof of Lemma 7.
The result is trivial if $C$ is planar, since any $\b{y}$ not lying in
the same plane as $C$ may be used.  We therefore assume that $C$ is
nonplanar. 
For the proof we shall find a plane $P$, given by an equation 
$\x.\b{a}=0$, so that $P$ intersects $C$ in exactly $d$ points
$\x_{i}$, say. We shall want $\b{a}$ to be a nonzero integer vector
satisfying $|\b{a}|\ll 1$.
We first demonstrate that this will suffice for our result.
We choose the vector $\b{y}$ to correspond to a point in the plane 
$P$, not on one of the lines 
$\langle\x_{i},\x_{j}\rangle$ for $i\not=j$.
Theorem 1 shows that this is possible with
$|\b{y}|\ll 1$. We then choose any integer vector $\b{c}$ for which
$\b{y}.\b{c}\not=0$ and $|\b{c}|\ll 1$, via a further application
of Theorem 1.  The projection $\pi$ from $C$ along $\b{y}$ onto the 
plane $\x.\b{c}=0$
is then a regular map, as $\b{y}$ is not on $C$, so that its image
$\pi(C)$ is
an irreducible curve.  Moreover the points
$\pi(\x_{1}),\ldots,\pi(\x_{d})$ are distinct, and lie on the
intersection of the curve $\pi(C)$ and the line $\pi(P)$. 
Thus $\pi(C)$ has degree at least $d$.
On the other hand, if $\pi(C)$ had degree greater than $d$ 
there would be a line $L$
intersecting $\pi(C)$ in more than $d$ points.  The inverse image
$\pi^{-1}(L)$ would then be a plane intersecting $C$ in more
than $d$~points, which is impossible, since $C$ is nonplanar.  
If the fibre over a point of $\pi(C)$
contained more than $d$~points, this would produce a line meeting $C$
in more than $d$~points.  Any plane containing this line would meet
$C$ in more than $d$ points, and hence would contain $C$.  This would
again contradict our assumption that $C$ is nonplanar.

The remainder of the proof is devoted to finding a suitable plane $P$.
According to B\'{e}zout's Theorem, in the form given by Harris [9;
Theorem 18.3], for example, it will suffice that $P$ passes
through none of the singular points of $C$, and is nowhere tangent to
$C$.

We begin by finding some equations of degree at most $d$, satisfied
on $C$.  We begin by choosing linearly independent vectors 
$\b{e}_{1},\ldots,\b{e}_{4}$ not lying on
$C$, and we change coordinates to use this as a new basis.  
The projection from $C$ along 
$\b{e}_{1}$ onto the plane spanned by $\b{e}_{2},\b{e}_{3},\b{e}_{4}$, 
is a regular
map, and the image is therefore an irreducible curve $C_{1}$, with
equation $f_{1}(x_{2},x_{3},x_{4})=0$.  The curve $C_{1}$ can have
degree at most $d$, by the argument above.  We therefore have an
absolutely irreducible equation $f_{1}(x_{2},x_{3},x_{4})=0$ of degree at most
$d$, satisfied everywhere
on $C$.  In the same way, we can produce absolutely irreducible
equations $f_{2}(x_{1},x_{3},x_{4})\break =0$ and
$f_{3}(x_{1},x_{2},x_{4})=0$, of degree at most $d$.
We shall think of each $f_{i}$
as being a form in $(x_{1},x_{2},x_{3},x_{4})$, being independent of
$x_{i}$. 

Let $I$ be the intersection
$$I:\;\;f_{1}(\x)=f_{2}(\x)=f_{3}(\x)=0.$$
If $I$ were to contain a component of dimension $2$, the polynomials
$f_{i}$, being absolutely irreducible, would have to be constant multiples of
each other.  This could only happen if they were each constant
multiples of $x_{4}$.  In this case however $C$ would be contained in
the plane $x_{4}=0$, contrary to assumption.
Now let $\Gamma$ be a component of $I$ of dimension $1$.  
We proceed to show that the $3\times 4$ 
matrix $M_{1}$, with rows $\nabla f_{1}(\x),\nabla f_{2}(\x)$ and 
$\nabla f_{3}(\x)$, has rank at least $2$ at a generic point $P_{0}$ of
$\Gamma$. Suppose, on the contrary, that $M_{1}$ has rank at most $1$ 
at $P_{0}$.
Since $\partial f_{i}/\partial x_{i}$ vanishes identically for $i=1,2$
and $3$, it then follows that there is some pair of indices $i\not=j$
for which
$$\frac{\partial f_{i}}{\partial x_{j}}(P_{0})=
\frac{\partial f_{j}}{\partial x_{i}}(P_{0})=0.$$
Suppose, to be specific, that $i=1,j=2$.  Then we have equations 
$$f_{1}(0,x_{2},x_{3},x_{4})=
\frac{\partial f_{1}}{\partial x_{2}}(0,x_{2},x_{3},x_{4})=0$$
and
$$f_{2}(x_{1},0,x_{3},x_{4})=
\frac{\partial f_{2}}{\partial x_{1}}(x_{1},0,x_{3},x_{4})=0,$$
holding on $\Gamma$.  If one of the partial derivatives, $\partial
f_{1}/\partial x_{2}$ say, vanishes identically, the form $f_{1}$
would take the shape $f_{1}(x_{3},x_{4})$.  Since $C$ lies on
$f_{1}=0$, it would follow that $C$ lies in a plane, which we assumed
was not the case.  We may therefore suppose that neither of the partial
derivatives above vanishes identically.
Since $f_{1}$ is absolutely irreducible, the first pair of
equations shows that $\Gamma$ must be a line through the point 
$(1,0,0,0)$.  Similarly the second pair of
equations shows that $\Gamma$ must be a line through the point 
$(0,1,0,0)$. Hence $\Gamma$
must be the line $x_{3}=x_{4}=0$.  The equation 
$f_{3}(x_{1},x_{2},0,x_{4})=0$ has to hold on this line, which implies that
$x_{4}|f_{3}(x_{1},x_{2},0,x_{4})$.  Since $f_{3}$ is irreducible,
this implies that $f_{3}(x_{1},x_{2},0,x_{4})=cx_{4}$.  However
$f_{3}=0$ is an equation for the original curve $C$, which was assumed
to be nonplanar.  This establishes our claim about the matrix $M_{1}$,
and shows that each one-dimensional component $\Gamma$ of $I$ 
contains only finitely many points where $M_{1}$ has rank at most $1$.
Indeed, since there are only finitely many components, and only
finitely many of these are points, we may conclude that there are only 
finitely many points $\x_{i}$ on $I$ which are either point components
of $I$ or for which $M_{1}$ has rank at most $1$.

Now let $\Delta_{i}(\b{a},\b{x})$ for $1\leq i\leq 16$ be the $3\times 3$
determinants formed from the matrix $M_{2}$ with rows $\nabla f_{1}(\x), 
\nabla f_{2}(\x),\nabla f_{3}(\x)$ and $\b{a}$, and consider the 
system of equations
\begin{equation}
\b{a}.\x=0,\;\;\;f_{i}(\x)=0\;\;\;(1\leq i\leq 3),\;\;\;
\Delta_{i}(\b{a},\b{x})=0\;\;\;(1\leq i\leq 16).\hskip.25in
\end{equation}
Choose a value of $\b{p}$ (not necessarily integral)
which does not lie on the dual variety
$\Gamma^{*}$, for any one-dimensional component $\Gamma$ of $I$, 
and such that $\b{p}.\x_{i}\not=0$ for each of points
$\x_{i}$ found above.  This is possible, since $\Gamma^{*}$ has dimension 
at most $2$ (see Harris [9; p.\  197]).  Then if $\x$ were a solution
to (4.1) with $\b{a}=\b{p}$, 
it would lie in the intersection $I$.  Moreover it cannot be one 
of the points $\x_{i}$, whence $\x$ lies on a curve $\Gamma$ in $I$, and 
$M_{1}$ has rank
at least $2$.  Each tangent space $T_{{\bf x}}(\Gamma)$ has projective 
dimension at
least one, and its elements are orthogonal to each of the 
$\nabla f_{i}(\x)$.
It follows that $M_{1}$ has rank exactly $2$, that $\Gamma$ is nonsingular at
$\x$, and that
$$T_{{\bf x}}(\Gamma)=\{\b{y}:\;\b{y}.\nabla f_{i}(\x)=0\;(1\leq i\leq 3)\}.$$
Since $\Delta_{i}(\b{p},\b{x})=0$ for each $i$ we see that $\b{p}$ is
in the linear span of the vectors $\nabla f_{i}(\x)$, whence
$$T_{{\bf x}}(\Gamma)\subset\{\b{y}:\;\b{y}.\b{p}=0\}.$$
This however contradicts our assumption that $\b{p}$ is not in
$\Gamma^{*}$.  Thus (4.1) has no
solutions $\x$ when $\b{a}=\b{p}$ .

Lemma 2 shows that there is a necessary and sufficient
condition for (4.1) to be solvable for $\x$, given by the vanishing of
a system of forms $G_{i}(\b{a})$.  These forms will have degrees which
are bounded in terms of $d$. The condition is nonempty, since (4.1) 
is not always solvable, as we
have shown.  We now make a linear change of variables to revert to our
original coordinate system.  Then, using Theorem 1, we can
find a nonzero integer vector $\b{a}\ll 1$ for which (4.1) has  
no solution.  Thus,
if $\x$ lies on $C$ and is also on the plane $\b{a}.\x=0$, we must
have $\Delta_{i}(\b{a},\b{x})\not=0$ for some $i$, since $f_{i}(\x)=0$
are amongst the equations for $C$.  It follows that $M_{2}$ has rank at
least $3$, and hence that $M_{1}$ has rank at
least $2$.  We can then deduce, as above, that $M_{1}$ has rank exactly
$2$, that $C$ is nonsingular at
$\x$, and that
$$T_{{\bf x}}(C)=\{\b{y}:\;\b{y}.\nabla f_{i}(\x)=0\;(1\leq i\leq 3)\}.$$
Since $M_{2}$ has strictly
larger rank than $M_{1}$, we see that $\b{a}$ is not
in the span of the vectors $\nabla f_{i}(\x)$, so that the tangent
space cannot be contained in the plane $\b{a}.\x=0$.  The plane
$\b{a}.\x=0$ therefore has the required properties, and Lemma~7 is
proved.

\section{Quadratic hypersurfaces}

This section is devoted to the proof of Theorem 2.  Our key tool is the
following result, for which see Heath-Brown [13; Theorem 3].

\specialnumber{8}\proclaim{Lemma}
Let $q$ be a nonsingular integral ternary quadratic form{\rm ,} 
with coefficients bounded
in modulus by $||q||,$ say{\rm .}  Suppose that the binary form 
$q(x_{1},x_{2},0)$ is also
nonsingular{\rm .}  Then for any integer $k$ the equation
$q(\b{x})=0$ has only
$O_{\ep}((||q||R)^{\ep})$ primitive integer solutions in the cube
$|x_{i}|\leq R,$ with $x_{3}=k.$
\endproclaim

We first prove a weaker version of Theorem 2, namely the estimate
\begin{equation}
N(B)\ll_{\ep} ||F||^{\ep}B^{n-2+\ep}.
\end{equation}
Having done this we shall use a technique similar to that developed
for Theorem~4, to deduce Theorem 2 itself.

To prove (5.1) we shall begin by making a change of variables,
$\x=M\b{y}$ to produce $F(M\b{y})=T(\b{y})$, say.  We shall write
$T_{ij}$ for the coefficients of $T$, so that the
$T_{ij}$ are quadratic polynomials in the entries $M_{ij}$.  We now
consider the function
$$f(M)=\det(M).T_{11}.\det(T_{ij})_{i,j\leq 2}.\det(T_{ij})_{i,j\leq 3}.$$
This does not vanish identically, since it is possible to choose $M$
so as to make $T$ diagonal, with at least $3$ nonzero entries.  Since
$f(M)$ is a form of degree $n+12$ in the entries $M_{ij}$ of the matrix
$M$, we see from Theorem 1 that there is
an integral matrix $M$, with $\max|M_{ij}|\ll 1$ such that
$f(M)\not=0$. If $\x\in\Z^{4}$ then $\det(M)\b{y}\in\Z^{4}$.  
Thus it suffices to consider solutions of $T(\b{y})=0$,
with $|\b{y}|\ll B$.  Here $||T||\ll ||F||$.

For any choice of $\b{u}=(y_{3},\ldots,y_{n})$ with $y_{i}\ll B$, 
we shall set
$$q(x,y,z)=T(x,y,z\b{u}).$$
The determinant of this form is a quadratic polynomial $D(\b{u})$, say.
Moreover, $D(\b{u})$ does not vanish identically, since 
$$D(1,0,0,\ldots,0)=\det(T_{ij})_{i,j\leq 3}\not=0,$$
by choice of $M$.  We also see that $q(x,y,0)$ is nonsingular,
because
$$\det(T_{ij})_{i,j\leq 2}\not=0,$$
again by choice of $M$.  Thus if $\b{u}$ is a value for which
$D(\b{u})\not=0$ then Lemma 8 shows that there are
$O_{\ep}((||F||B)^{\ep})$ possible values of $y_{1},y_{2}$ making
$T(\b{y})=0$.  This produces $O_{\ep}(||F||^{\ep}B^{n-2+\ep})$
solutions in total.  On the other hand, since $D(\b{u})$ does not
vanish identically, there can be only $O(B^{n-3})$ values of $\b{u}$
for which $D(\b{u})=0$, by Theorem 1.  
For each of these we can specify $y_{2}$ in
$O(B)$ ways, and then there are at most $2$ corresponding values of
$y_{1}$, since $T_{11}\not=0$, by choice of $M$.  There are therefore
$O(B^{n-2})$ solutions for which $D(\b{u})=0$, which completes the
proof of (5.1).

To derive Theorem 2 from (5.1) we shall adapt
the treatment of Theorem~4.  Let $\x^{(1)},\ldots,\x^{(N)}\in Z_{n}$
be the complete set of solutions of $F(\x)=0$ in
the region $|\x^{(i)}|\ll B$.  Set $M=n(n+1)/2$ for convenience, and 
consider the
$N\times M$ matrix $C$, whose $i^{\rm th}$ row consists of
the $M$ possible monomials of degree $2$ in the variables 
$x^{(i)}_{1},\ldots,x^{(i)}_{n}$.  Then if the vector 
$\b{f}\in\Z^{M}$ has entries which are the corresponding coefficients 
of $F$, we will have $C\b{f}=\b{0}$.  Since $\b{f}\not=\b{0}$ it 
follows that $C$ has rank at most $M-1$.  Thus $C\b{g}=\b{0}$ has a
nonzero integer solution~$\b{g}$, constructed out of the
sub-determinants of $C$.  It follows that there is such a $\b{g}$ with
$|\b{g}|\ll_{d} B^{2M-2}$. Let $G(\x)$ be the quadratic form 
corresponding to the vector $\b{g}$.  Then $G(\x)$ and $F(\x)$
have $N$ common zeros, namely the vectors $\x^{(i)}$.  If $G(\x)$ is a
rational multiple of $F(\x)$ then
$$N(F;B)\leq N(G;B)\ll_{\ep} ||G||^{\ep}B^{n-2+\ep},$$
by (5.1).  In this case we have $N(F;B)\ll_{\ep} B^{2+\ep}$, as
required, on re-defining $\ep$.

If $G(\x)$ is not a rational multiple of $F(\x)$ then the points
$\x^{(i)}$ satisfy $F(\x)=G(\x)=0$.  As above, we may apply a linear
transformation so that $F$ contains the term $x_{1}^{2}$ with nonzero
coefficient.  We can then eliminate $x_{1}$ from the equations 
$F(\x)=G(\x)=0$ to deduce that
$H(x_{2},\ldots,x_{n})=0$, for some nonzero form $H$ of degree at
most $4$.  Theorem 1 shows that this has $O(B^{n-2})$ solutions in the
relevant region, and for each of these solutions $(x_{2},\ldots,x_{n})$
the equation $F(\x)=0$
determines at most two values of $x_{1}$.  It follows that $N\ll
B^{n-2}$ in this case, and Theorem 2 follows. \pagebreak

\section{General surfaces}

In this section we shall consider Theorems 6, 7 and 9. We begin with
Theorem 6. 
We shall apply part (iv) of Lemma 1 in the case in which $n=4$, so that 

\begin{equation}
\b{y}\ll B^{1/3}.  
\end{equation}
The points on $F(\x)=0$ which also 
lie in the plane $\x.\b{y}=0$ are in one-to-one correspondence with points 
on a curve 
$G_{\bf y}(\lambda_{1},\lambda_{2},\lambda_{3})=0,$ where
$$G_{\bf y}(\lambda_{1},\lambda_{2},\lambda_{3})=
F\lrp{\lambda_{1}\x^{(1)}+\lambda_{2}\x^{(2)}+\lambda_{3}\x^{(3)}}.$$
Moreover primitive points on $F=0$ correspond to primitive points on
$G_{\bf y}=0,$ and vice-versa.  If $\x=\sum\lambda_{j}\x^{(j)}$ lies
in the box $\max|x_{i}|\leq B$, then Lemma 1, part (iii), 
yields $|\lambda_{j}|\ll
B|\x^{(j)}|^{-1}$.  We then apply Theorem 3 with
$B_{j}=cB|\x^{(j)}|^{-1}$, for a suitable constant $c>0$. 
If we order the indices so that $|\x^{(1)}|\geq
|\x^{(2)}|\geq |\x^{(3)}|$, we will
have $T\gg B^{d}|\x^{(1)}|^{-d}$,
whence 
\begin{equation}
N\lrp{G_{\bf y};c\frac{B}{|\x^{(1)}|},c\frac{B}{|\x^{(2)}|},
c\frac{B}{|\x^{(3)}|}}\ll_{c,\ep} B^{2/d+\ep}
(|\x^{(2)}|.|\x^{(3)}|)^{-1/d},\hskip.4in
\end{equation}
providing that $G_{\bf y}$
is irreducible over $\Q$.  

We now sum over the possible vectors $\b{y}$,
counting them according to the values of the various $\x^{(j)}$.  
Consider the case
in which $C_{j}<|\x^{(j)}|\leq 2C_{j}$ for $1\leq j\leq 3$.  The vector
$\b{y}$ lies in the integer lattice defined by $\b{y}.\x^{(3)}=0$, and
this lattice has determinant $|\x^{(3)}|$, by Lemma 1, part (i).  In our
situation we take $\b{y}\ll Y=C_{1}C_{2}C_{3}$, in view of (2.2),
whence part (v) of Lemma 1 shows that the number of possible vectors 
$\b{y}$ is
$O(Y^{3}|\x^{(3)}|^{-1})$.  Thus there
are $O(C_{1}^{3}C_{2}^{3}C_{3}^{2})$ possible values of $\b{y}$ for
each $\x^{(3)}$.  We sum this over the $O(C_{3}^{4})$ possible vectors
$\x^{(3)}$, and conclude that there are $O(C_{1}^{3}C_{2}^{3}C_{3}^{6})$ 
values of $\b{y}$ for which $C_{j}<|\x^{(j)}|\leq 2C_{j}$ for 
$1\leq j\leq 3$.  For each such $\b{y}$ we have
$$N\lrp{G_{\bf y};c\frac{B}{|\x^{(1)}|},c\frac{B}{|\x^{(2)}|},
c\frac{B}{|\x^{(3)}|}}\ll_{c,\ep} 
B^{2/d+\ep}(C_{2}C_{3})^{-1/d},$$
by (6.2), producing a total contribution
$$\ll_{\ep} B^{2/d+\ep}C_{1}^{3}C_{2}^{3-1/d}C_{3}^{6-1/d}.$$
Since the indices are ordered with $C_{1}\gg C_{2}\gg C_{3}$, and
$$C_{1}C_{2}C_{3}\ll |\x^{(1)}|.|\x^{(2)}|.|\x^{(3)}|\ll |\b{y}|\ll B^{1/3},$$
by (2.2) and (6.1), we obtain an estimate
\begin{eqnarray*}
&\ll_{\ep} & B^{2/d+\ep}(C_{1}C_{2}C_{3})^{4-2/3d}\\
&\ll_{\ep} & B^{2/d+\ep}B^{4/3-2/9d}\\
&\ll_{\ep} & B^{4/3+16/9d+\ep},
\end{eqnarray*}
for the contribution to $N(B)$ corresponding to the case in which $G_{\bf y}$
is irreducible over $\Q$, and $C_{j}<|\x^{(j)}|\leq 2C_{j}$.  
Finally we let the
$C_{j}$ run over powers of $2$ and sum the resulting bounds to obtain
an estimate which we state formally as follows.

\specialnumber{9}\proclaim{Lemma}
The contribution to $N(B)$ corresponding to those vectors $\b{y}$ for 
which $G_{\bf y}$ is irreducible over $\Q$ is 
$O_{\ep}(B^{4/3+16/9d+\ep})${\rm .}
\endproclaim

We must now tackle the case in which $G_{\bf y}$ is reducible over $\Q$.  If\break
$G_{\bf y}(\lambda_{1},\lambda_{2},\lambda_{3})=0,$ then we must have
$H(\lambda_{1},\lambda_{2},\lambda_{3})=0$ for some factor $H$ of
$G_{\bf y}$. We may suppose that $H$ is irreducible over $\Q$, though
not necessarily absolutely irreducible.
Of course, any solution corresponding to a linear factor $H$ produces
a point $\x$ lying on a line in the surface $F=0$ which is defined over~$\Q$.
We next dispose of the case in which $H$ has degree $d'\geq 3$. 
Here the analysis leading up to Lemma 9 goes through just as before,
and leads to a contribution 
\begin{equation}
\ll_{\ep}B^{4/3+16/9d'+\ep}
\ll_{\ep}B^{52/27+\ep}.
\end{equation}
We turn now to the case in which there is a quadratic factor.
We shall assume in what follows that $d\geq 3$.  Lemma 3 shows that
there is a set of conditions $E_{m}(\b{y})=0$ which are necessary and
sufficient for $G_{\bf y}$ to have a quadratic factor.
In general an elimination
procedure of the above type may lead to an empty set of equations
$E_{m}=0$.  However in our case this does not happen, since the 
generic plane section of the
surface $F=0$ is known to be irreducible, see Harris 
[9; Proposition 18.10].  
At least one of the forms $E_{m}$ must therefore be nonzero, and 
we may therefore conclude as follows.

\specialnumber{10}\proclaim{Lemma}
There is a nonzero form $E(\b{y})$ with degree bounded in terms of
$d${\rm ,} 
such that if $G_{\bf y}$ has a quadratic factor{\rm ,} then $E(\b{y})=0${\rm .}
\endproclaim

It should be stressed that the only respect in which this differs from
the statement that the generic plane section of the
surface $F=0$ is irreducible, lies in our control over the degree of
$E$.

We can now apply Theorem 1 to show that there are
$O(Y^{3})$ vectors $\b{y}$ with $Y<|\b{y}|\leq 2Y$, such that $G_{\bf y}$ has
a quadratic factor, $H$ say.  Then if $H$ were singular, but
irreducible over $\Q$, we would find that
$H(\lambda_{1},\lambda_{2},\lambda_{3})=0$ has $O(1)$ primitive solutions 
$(\lambda_{1},\lambda_{2},\lambda_{3})$.  If $H$ is nonsingular, 
we may apply Theorem 3, with $B_{i}\ll B|\x^{(i)}|^{-1}$, to deduce
that
\begin{eqnarray*}
N\lrp{H;c\frac{B}{|\x^{(1)}|},c\frac{B}{|\x^{(2)}|},c\frac{B}{|\x^{(3)}|}}
&\ll_{c,\ep} &B^{1+\ep}(|\x_{1}|.|\x_{2}|.|\x_{3}|)^{-1/3}\\
&\ll_{c,\ep} &B^{1+\ep}|\b{y}|^{-1/3},
\end{eqnarray*}
by (2.2).  The range \pagebreak $Y<|\b{y}|\leq 2Y$ therefore contributes 
$O_{\ep}(B^{1+\ep}Y^{8/3})$.  Thus if we sum $Y$ over
powers of 2, with $Y\ll B^{1/3}$, we obtain a contribution
$O_{\ep}(B^{17/9+\ep})$.  If we combine this with the
bounds given by (6.3) and by Lemma~9, we obtain the assertion of Theorem
6. 

We turn now to Theorem 7. Here our starting point is the case $n=4$ of
Theorem 14, which shows that every point $\x$ on the surface $F(\x)=0$
which lies in the cube $|\x|\leq B$, must also satisfy one of the
equations $F_{j}(\x)=0$.  Here 
$$j\leq k\ll_{\ep} B^{3/\sqrt{d}+\ep}\log^{5}||F||.$$
The intersection $F(\x)=F_{j}(\x)=0$ consists of at most $dD$ curves $C$,
with degrees at most $dD$.  If $C$ is a line, not defined over $\Q$,
it contains at most one rational point.  We can therefore suppose that
$C$ has degree at least $2$. To estimate the number of points on such a
curve $C$, we apply Theorem 5.  Thus each curve contributes 
$O_{d,\ep}(B^{1+\ep})$ points,
so that
$$N_{1}(F;B)\ll_{\ep}B^{1+3/\sqrt{d}+\ep}\log^{5}||F||.$$

We proceed to show that the factor $\log^{5}||F||$ can be eliminated from
this estimate, by the method used in Section 5 for proving the case $d=2$ 
of Theorem 9.  We take $\x^{(1)},\ldots,\x^{(N)}$ to be the
complete set of solutions of $F(\x)=0$ in the region $|\x^{(i)}|\ll
B$, excepting any that lie on lines in the surface which are defined
over $\Q$.  Proceeding as
before, we reach two possible cases.  In the first case, when the form
$G$ is a constant multiple of $F$, we deduce that
\begin{eqnarray*}
N_{1}(F;B)&=&N_{1}(G;B)\\
&\ll_{\ep}&B^{1+3/\sqrt{d}+\ep}\log^{5}||G||\\
&\ll_{\ep}&B^{1+3/\sqrt{d}+\ep}\log^{5}B,
\end{eqnarray*}
which suffices for Theorem 7.  In the second case, all the points
$\x^{(1)},\ldots,\x^{(N)}$ lie on one of at most $O_{d}(1)$ curves $C$
of degree at most $d^{2}$, lying in the surface.  By definition of
$N_{1}(F;B)$ these curves have degrees $\delta\geq 2$.  Thus Theorem 5
shows that each curve contains at most
$O_{\delta,\ep}(B^{1+\ep})$ points, whence
$$N=N_{1}(F;B)\ll_{\ep,d}B^{1+\ep}$$ 
in this case.  This completes the proof of Theorem 7.

Turning finally to Theorem 9, we see from Theorem 2 that it suffices
to take
$d\geq 3$.  In view of Theorem 6, we will have to estimate the contribution
from  lines on the surface $S$ given by $F=0$. These lines 
correspond to points in the
Grassmannian $\Grass(1,3)=G$, say.  Indeed those lines that lie in the 
surface
$S$ correspond to points of an algebraic subset $V$, say of
$G$, (the  Fano variety $F_{1}(S)$, see Harris [9; Example 6.19]).  
The set $V$ is
defined by $O_{d}(1)$ equations of degree at most $d$.  The lines that
lie in a plane $P$ correspond to points \pagebreak  on a plane
$P'\subset G$.  For a generic
plane $P\subset\Proj^{3}$ the intersection $P\cap S$ is irreducible,
(see Harris [9; Proposition 18.10]) and so contains no lines.  Hence
there is a plane $P\subset\Proj^{3}$ for which the corresponding $P'$ 
is disjoint from $V$.

If we choose coordinates so that $P$ consists of points
$(0,x,y,z)$, then the Pl\"{u}cker coordinates $p_{ij}$ of the lines in
$P$ all have $p_{12}=p_{13}=p_{14}=0$.  We now choose $A,B$ such that 
$F(0,A,B,1)\not=0$.  This is clearly possible since 
$x_{1}\nmid F(x_{1},x_{2},x_{3},x_{4})$.  Then the intersection of
$G$, given by
$$p_{12}p_{34}-p_{13}p_{24}+p_{14}p_{23}=0,$$
with the linear space $L$, given by $p_{12}=Ap_{14}$ and
$p_{13}=Bp_{14}$, will be the union of the plane $P'$ with a second
plane $P''$, given by the equations $p_{12}=Ap_{14}$, $p_{13}=Bp_{14}$
and $Ap_{34}-Bp_{24}+p_{23}=0$.  This second plane corresponds to the
set of lines passing through the point $(0,A,B,1)$.  By construction
none of these lines lie in $S$, so that $P''$ is also disjoint from $V$.

It follows that  
$$L\cap V=L\cap G\cap V=(P'\cup P'')\cap V=\emptyset,$$
from which we conclude that every component of $V$ has dimension at
most $1$.  A line in $G$ corresponds to the set of lines in $\Proj^{3}$
which lie in a given plane and pass through a given point.  At most
finitely many of these can be contained in $S$, whence $V$ cannot
contain a line.  

Any line that does not pass through two distinct rational points can
contribute at most $1$ to $N(F;B)$.  There are $O(P^{4/3})$ lines to
consider so the total contribution from lines which are not defined
over $\Q$ is
$O(P^{4/3})$. We therefore focus our attention on those lines which
are defined over $\Q$.  These
correspond to rational points lying in
$V$.  We re-scale these so as to be primitive integral points, and count
the number of such points with height at most $Y$, say.  To do this we shall
investigate projections of $V$ onto
various linear spaces.  Choose a vector $\b{p}\in\Proj^{5}$ not lying
on $V$, and a hyperplane $H$ not containing $\b{p}$.  Let $\pi$ be the 
projection from $V$ to $H$ along $\b{p}$.  Then $\pi$ is a regular 
map, and the image
$\pi(V)$ is therefore a closed algebraic set, with components of
dimension at most $1$.  One can produce a set
of defining equations for $\pi(V)$ via elimination theory, and one sees
that there will be $O_{d}(1)$ equations, with degrees bounded in terms
of $d$.  To be specific, let $f_{i}(\x)=0$ be a suitable set of
defining equations for $V$, and let $\b{h}\in H$.  
According to Lemma 2, the system of
equations $f_{i}(\lambda\b{p}+\mu\b{h})=0$ will have a nonzero solution
$\lambda,\mu$, if and only if $\b{h}$ satisfies a system of polynomial
equations $E_{j}(\b{h})=0$.  Since $\b{p}\not\in V$, we must have
$\mu\not=0$ and $\b{h}\not=0$ in any such solution.  The projection
$\pi(V)$ is therefore given by the equations $E_{j}(\b{h})=0$.  If
these are not homogeneous, then $\b{h}$ must clearly be a zero of each
of their homogeneous components.  One therefore obtains
in this way a collection of $O_{d}(1)$ conditions, of degrees
$O_{d}(1)$. 

If $C$ is an irreducible component of $V$, then $\pi(C)$ will
be an irreducible component of $\pi(V)$.  Since $C$ cannot be a line,
it follows that $\pi(C)$ cannot be a point.  Moreover, if $\pi(C)$ is
a line, then $C$ is planar, lying in a plane $P_{C}$, say, containing 
$\b{p}$. We now choose $\b{p}$ not lying in any of the planes $P_{C}$,
nor on $V$.  Thus it suffices for some form of degree $O_{d}(1)$ to be
nonvanishing at $\b{p}$.  In view of Theorem 1 we can choose
an integral point of this type, such that $|\b{p}|\ll 1$.  
Similarly we can choose a
hyperplane $H$ given by $\b{c}.\b{h}=0$ with $\b{c}$ integral,
so that $|\b{c}|\ll 1$ and
$\b{p}\not\in H$.  It then follows that every component $C$ of $V$
projects to a curve in $H$ which is not a line.  Moreover we can choose
coordinates in $H$ so that points of height at most $Y$ in $\Proj^{5}$
project to points of height $O(Y)$.  Since no component $C$
projects to a point, it follows that the inverse image of any point on
$\pi(C)$ contains $O(1)$ points.

In order to estimate the number of points in $V$ it therefore suffices
to estimate the number of points in $\pi(V)$.  Clearly we may iterate
this process, reducing the problem to one of points on a plane curve.
In this case Theorem 3 gives a bound $O_{\ep}(Y^{1+\ep})$, so that
we may conclude that $V$ itself contains $O_{\ep}(Y^{1+\ep})$ points
of height at most $Y$.

Each line $L\subset S$ which is defined over $\Q$ 
intersects $\Z^{4}$ in a lattice $\sLambda$,
say, of rank 2.  If the lattice has determinant $\Delta$, then
Lemma 1, part (vi), shows that the
line will contain $O(1+B^{2}/\Delta)$ points of $Z_{4}$ from the cube\break
$\max |x_{i}|\leq B$.  However the determinant $\Delta$ is merely the
height of the corresponding Pl\"{u}cker coordinate vector, which we
take to be primitive.  The lines in which we are interested arise from
the intersection of the surface $S$ with various planes  $\x.\b{y}=0$,
with $|\b{y}|\ll B^{1/3}$.  There are therefore $O(B^{4/3})$ such
lines, so that lines with $\Delta\geq B^{2}$ contribute $O(B^{4/3})$ 
to $N(F;B)$. Moreover, as we have just shown, there are
$O_{\ep}(Y^{1+\ep})$ lines with $Y<\Delta\leq 2Y$.  When $Y\ll B^{2}$,
such lines therefore contribute $O_{\ep}(B^{2}Y^{\ep})$.
Finally we may sum over values of $Y$ running over powers of $2$, to
obtain an overall contribution $O_{\ep}(B^{2+2\ep})$ from lines in
$S$.  This completes the proof of Theorem 9.

We conclude by remarking that the above method fails for $d=2$ only
because the analogue of Lemma 9 would contain an exponent
$4/3+16/9d=20/9>2$. The treatment of points on lines in the surface
still applies satisfactorily.

\section{Binary forms}

This section is devoted to the proof of Theorem 8.  It will be
convenient to make a
linear change of variable so that $G(1,0)\not=0$.  Clearly this has no
effect on the conclusion of Theorem 8.  Throughout this section, all
implied constants may depend on the form $G$.  We shall not mention
this dependence explicitly. We begin by defining
\begin{eqnarray*}
   S(X,C)&  =&\#\{(x,y)\in\Z^{2}: 1\leq G(x,y)\leq X,\\
&&\hskip.5in C<\max(|x|,|y|)\leq 2C,
\;{\rm h.c.f.}(x,y)=1\},
\end{eqnarray*}
subject to the assumption that $C\gg X^{1/d}$.
We observe that if $x,y$ is counted by $S(X,C)$ then there is some
factor $x-ay$ of $G(x,y)$ such that $|x-ay|\ll X^{1/d}$.  Thus if
$C\geq cX^{1/d}$ with a sufficiently large constant $c$, we will have 
$$C\ll |x-a'y|\ll C$$
for every factor with $a'\not=a$.  If $x-ay$
divides $G$ with multiplicity $e$ it then follows that 
$|x-ay|^{e}C^{d-e}\ll X$.
If $a$ is irrational we have $|x-ay|\gg_{\ep} C^{-1-\ep}$, by Roth's theorem. If $a$  is
rational, we cannot have $x-ay=0$,
since $n\not=0$. It follows that $|x-ay|\gg 1$ if $a$ is rational.  
In either case we may use our
assumption that $e\leq (d-1)/2$ to deduce that
$$C^{1-d\ep}\leq C^{d-2e-e\ep}\ll_{\ep}|x-ay|^{e}C^{d-e}\ll X,$$
whence $C\ll X^{2}$, say.  Thus $S(X,C)=0$ unless $C\ll X^{2}$.

We proceed to estimate the contribution to $S(X,C)$ corresponding to a
particular value of $a$.  Such a contribution arises from 
primitive lattice points in the
parallelogram $|y|\leq 2C,\;|x-ay|\ll (XC^{e-d})^{1/e}$.  According
to Lemma~1, part (vii), there are 
$$\ll 1+C(XC^{e-d})^{1/e}=1+C^{2}(XC^{-d})^{1/e}\leq
1+C^{2}(XC^{-d})^{2/(d-1)}$$ 
such points, using once more the assumption that $e\leq (d-1)/2$.  
It therefore follows that
$$S(X,C)\ll 1+(X/C)^{2/(d-1)}.$$
We may now sum up for values of
$C\ll X^{2}$, running over powers of $2$, to conclude that if
\begin{eqnarray*}
&&{S'(X,C)}\\
&&\quad=\#\{(x,y)\in\Z^{2}: 1\leq G(x,y)\leq X,\;\max(|x|,|y|)> C,\;
{\rm h.c.f.}(x,y)=1\}
\end{eqnarray*}
then
$$S'(X,C)\ll\log X+(X/C)^{2/(d-1)},$$
for $C\gg X^{1/d}$.

We now set
$$r(n)=\#\{(x,y)\in\Z^{2}: n=G(x,y)\}$$
and
$$r_{1}(n;C)=\#\{(x,y)\in\Z^{2}: n=G(x,y),\;\max(|x|,|y|)\leq C\},$$
$$r_{2}(n;C)=\#\{(x,y)\in\Z^{2}: n=G(x,y),\;\max(|x|,|y|)> C\},$$
where $C\gg X^{1/d}$.  Then
\begin{eqnarray}
\sum_{n\leq X}r_{2}(n;C)&=&\sum_{h\ll X^{1/d}}
S'\lrp{\frac{X}{h^{d}},\frac{C}{h}}\\
&\ll&\sum_{h\ll X^{1/d}}\left\{\log X+
\lrp{\frac{X}{h^{d}}\frac{h}{C}}^{2/(d-1)}\right\}\nonumber\\
&\ll& X^{1/d}\log X+\lrp{\frac{X}{C}}^{2/(d-1)}\sum_{h\ll X^{1/d}}h^{-2}
\nonumber\\
&\ll& X^{1/d}\log X+\lrp{\frac{X}{C}}^{2/(d-1)}.\nonumber
\end{eqnarray}
The sum $\sum r_{1}(n;C)$ is trivially $O(C^{2})$, whence
$$\sum_{n\leq X}r(n)\ll C^{2}+X^{1/d}\log X +\lrp{\frac{X}{C}}^{2/(d-1)}
\ll X^{2/d},$$
on choosing $C=cX^{1/d}$ with an appropriate constant $c$.  This bound
shows that there are $O(X^{2/d})$ positive integers $n\leq X$
represented by $G$.

If $n$ has two inequivalent representations $G(x,y)=n$, then either
$r_{2}(n;C)$ is positive, or there is a point
$(x_{1},x_{2},x_{3},x_{4})$ on the
surface 
$$E(\x)=G(x_{1},x_{2})-G(x_{3},x_{4})=0,$$
satisfying
$|x_{i}|\leq C$, but for which $(x_{1},x_{2})$ and $(x_{3},x_{4})$ are
not related by an automorphism.  We shall let $\cl{N}(C)$ denote the
number of such points.  We now claim that
\begin{equation}
\cl{N}(C)\ll_{\ep}C^{52/27+\ep},
\end{equation}
and that the form $G$ has $O(1)$ automorphisms.  It will then
follow that
$$\sum_{n\leq X}r_{1}(n;C)^{2}\leq \cl{N}(C)+O\lrp{\sum_{n\leq X}r_{1}(n;C)}
\ll C^{2}.$$
If $G(1,0)>0$ we trivially have
$$\sum_{n\leq X}r_{1}(n;C)\gg C^{2}$$
if $C=cX^{1/d}$ with a sufficiently small constant $c$, since then
$\max(|x|,|y|)\leq C$ implies $|G(x,y)|\leq X$, and
a positive proportion of such pairs $x,y$ will have $G(x,y)>0$.  It
now follows via Cauchy's inequality that $r_{1}(n;C)>0$ for $\gg
C^{2}$ positive integers $n\leq X$.  Thus the number of such integers
represented by $G$ has exact order $X^{2/d}$, as claimed in Theorem 8.

For integers $n$ with two or more
essentially different representations, we observe as above that either 
$r_{2}(n;C)>0$ or the representations are counted by $\cl{N}(C)$.
Thus the number of such integers will be
\begin{eqnarray*}
&\leq& \cl{N}(C)+\sum_{n\leq X}r_{2}(n;C)\\
&\ll_{\ep}&C^{52/27+\ep}+X^{1/d}\log X +\lrp{\frac{X}{C}}^{2/(d-1)}\\
&\ll_{\ep}&X^{52/(1+26d)+\ep},
\end{eqnarray*}
by (7.1) and (7.2), on choosing $C=X^{27/(1+26d)}$.  This completes
the proof of Theorem 8, subject to the claims made above.

We therefore turn to the consideration of integral points on the surface
$E(\x)=G(x_{1},x_{2})-G(x_{3},x_{4})=0$, in the cube $\max|x_{i}|\leq C$.  
We shall show that $E$ has no rational linear or quadratic factor.
Suppose to the contrary that there is such a factor.
Set $x_{1}=x,\,x_{2}=0,\,
x_{3}=a$ and $x_{4}=1$.  Then $E(\x)$ reduces to $Ax^{d}-B$, where 
$A=G(1,0)$ and
$B=G(a,1)$  Thus $Ax^{d}-B$ has a rational linear, or quadratic 
factor, for every
integral $a$.  It follows that $B/A$ is always an exact $d^{\rm th}$ power,
or, if $d$ is even, an exact $d/2^{\rm th}$ power.  This implies that $G(x,y)$ is
a perfect $d^{\rm th}$ power, or, if $d$ is even, a perfect $d/2^{\rm th}$ power.
However our assumption about the multiplicity of the factors of $G$
shows that this is impossible.

We can now apply Theorem 6 to each factor of $E(\x)$ to show that
$$N_{1}(E;Y)\ll_{\ep}Y^{52/27+\ep}.$$
If $\cl{N}^{(*)}(C)$ denotes the
number of integral zeros of $E$, not necessarily primitive, lying in
the cube $|x_{i}|\leq C$, but not on any line in the surface $E=0$,
then we conclude that
\begin{eqnarray*}
\cl{N}^{(*)}(C)&=&1+\sum_{h\ll C}N_{1}(E;B/h)\\
&\ll_{\ep}&1+\sum_{h}(C/h)^{52/27+\ep}\\
&\ll_{\ep}&C^{52/27+\ep}.
\end{eqnarray*}
The contribution to $\cl{N}(C)$ from points not lying
on lines in the surface $E=0$ is thus $O_{\ep}(C^{52/27+\ep})$, in
accordance with (7.2).

Lines in $\Proj^{3}$ may be classified into two types, given
respectively by pairs of equations
$$a_{1}x_{1}+a_{2}x_{2}+a_{3}x_{3}+a_{4}x_{4}=0,\;\;
b_{3}x_{3}+b_{4}x_{4}=0,$$
and
$$x_{1}=a_{1}x_{3}+a_{2}x_{4},\;\;x_{2}=a_{3}x_{3}+a_{4}x_{4}.$$
A little thought shows that if a line of the first type lies in the 
surface $E(\x)=0$, then the points on it must satisfy
$G(x_{1},x_{2})=G(x_{3},x_{4})=0$, since $G$ is not a $d^{\rm th}$ power.
Such points therefore correspond to the
excluded value $n=0$. Similarly, lines of the 
second type for which $a_{1}a_{4}=a_{2}a_{3}$ also produce values with
$G(x_{1},x_{2})=G(x_{3},x_{4})=0$.  The remaining
lines produce automorphisms 
\begin{equation}
G(a_{1}x+a_{2}y,a_{3}x+a_{4}y)=G(x,y).
\end{equation}
Indeed, if the $a_{i}$ are rational, the corresponding points produce
equivalent solutions of $G(x,y)=n$, in the sense of Theorem 8.
Thus points counted by $\cl{N}(C)$ which lie on lines in the surface
$E=0$ must lie on lines that correspond to irrational automorphisms.

We claim that there are only finitely many automorphisms,
rational or irrational.  Since 
any line which is not defined over $\Q$ contains at most $O(C)$ integral 
points (not necessarily primitive), we will be able to conclude 
that lines
corresponding to irrational automorphisms contribute $O(C)$ to $\cl{N}(C)$.
We will then have $\cl{N}(C)\ll_{\ep}C^{52/27+\ep}$, as required for
(7.2). 

It remains to prove that there are finitely many automorphisms.
The automorphisms of $G$ form a group, which acts on the roots of the
polynomial $G(x,1)$.  Specifically, the automorphism (7.3) maps a root
$\alpha$ by
\begin{equation}
\alpha\mapsto\frac{a_{1}\alpha+a_{2}}{a_{3}\alpha+a_{4}}.  
\end{equation}
The condition
on the multiplicity of the factors of $G$ implies that there are at
least three different roots $\alpha$.  For an automorphism
which fixed every root $\alpha$, one would have a quadratic equation
$$a_{3}\alpha^{2}+(a_{4}-a_{1})\alpha-a_{2}=0$$
with three distinct roots.  This would entail $a_{1}=a_{4}$ and
$a_{2}=a_{3}=0$.  One would then deduce from (7.3) that the common value
of $a_{1}$ and $a_{4}$ must be a $d^{\rm th}$ root of unity.  If we factor
the group of automorphisms by the subgroup consisting of scalar multiples
of $d^{\rm th}$ root of unity, the quotient still acts on the roots $\alpha$
by the formula (7.4), and the action is now faithful.  The quotient
group is thus isomorphic to a subgroup of the symmetric group
$S_{d}$.  We conclude that there are at most $d.d!$ automorphisms.

\section{Nonsingular surfaces}

In this section we shall prove Theorems 10 and 11.  The argument 
for Theorem 10 begins in
exactly the same way as for Theorem 6, and indeed, Lemma 9 shows that
the contribution from planes $\x.\b{y}=0$ for which $G_{\bf y}$ is
irreducible over $\Q$ is satisfactory.

We therefore consider the possibility that $G_{\bf y}$ factors.  In
this case $G_{\bf y}$ is a singular form, so that $\b{y}$ lies on the
dual surface $\hat{F}(\b{y})=0$. We proceed to show that $\hat{F}$
cannot be linear.  Since $\hat{F}(\nabla F(\x))$
vanishes on $F(\x)=0$ we have $F(\x)|\hat{F}(\nabla F(\x))$.  Hence if
$\hat{F}$ were linear we would deduce that $\hat{F}(\nabla F(\x))$
must vanish identically.  Taking $\hat{F}(\b{y})$ to have the shape
$\b{h}.\b{y}=0$ we would then have $\b{h}.\nabla F(\x)=0$ identically in
$\x$.  On taking the partial derivative with respect to $x_{j}$, say,
we conclude that
$$\sum_{i=1}^{4}h_{i}
\frac{\partial^{2}F(\x)}{\partial x_{i}\partial x_{j}}=0.$$
If we substitute $\b{h}$ for $\x$ this yields
$$\frac{\partial F}{\partial x_{j}}(\b{h})=0,$$
and since $j$ is arbitrary we have $\nabla F(\b{h})=0$.  This would
contradict the assumption that $F$ is nonsingular, so that $\hat{F}$
cannot be linear.

We may now apply Theorem 9, to show that 
there can be $O_{\ep}(B^{2/3+\ep})$ 
possible vectors $\b{y}$ with $|\b{y}|\ll B^{1/3}$.  Let $H$ be a
factor of $G_{\bf y}$ irreducible over $\Q$, and suppose that $H$ has degree
$e$.  According to Theorem 3 and Corollary 1,
the number of points $(\lambda_{1},\lambda_{2},\lambda_{3})$,
with $|\lambda_{i}|\ll B$, which satisfy
$H(\lambda_{1},\lambda_{2},\lambda_{3})=0$, will be $O_{\ep}(B^{2/e+\ep})$.
When $e\geq 3$, such factors $H$ produce a total contribution
$O_{\ep}(B^{4/3+\ep})$, after allowing for
$O_{\ep}(B^{2/3+\ep})$  possible $\b{y}$.  This is satisfactory.  For
$e=2$ we get an estimate $O_{\ep}(B^{5/3+\ep})$ in an analogous
fashion.  This too will be satisfactory providing that $d\leq 5$.

It remains to consider the possibility of quadratic factors $H$ of
$G_{\bf y}$, when $F$ has degree $d\geq 6$.  
However Theorem 12 then shows that the surface $F=0$ contains
$O(1)$ plane quadrics, so that $G_{\bf y}$ can have a quadratic
factor for at most $O(1)$ values of $\b{y}$.  Each such factor
produces $O_{\ep}(B^{1+\ep})$ points, giving a total contribution to
$N(F;B)$ of $O_{\ep}(B^{1+\ep})$.  This suffices for Theorem 10.

We turn now to Theorem 11. Our argument begins in precisely the same
way as was used in Section 6, for Theorem 7.  Thus every point on the
surface $F(\x)=0$, contained in the cube $|x_{i}|\leq B$, lies on one
of
$$\ll_{\ep} B^{3/\sqrt{d}+\ep}\log^{5}||F||$$
curves $C$.  Moreover the degree $\delta$ of any such curve is $O(1)$.
To estimate the number of points on such a
curve $C$, we again apply Theorem 5, to conclude that there are 
$O_{\ep}(B^{2/(d-1)+\ep})$
points lying on $C$ whenever $\delta\geq d-1$. Thus there are
a total of  $O_{\ep}(B^{3/\sqrt{d}+2/(d-1)+\ep}\log^{5}||F||)$ points 
lying on the available collection of curves $C$ of degree $d-1$ or more.
This shows that
\begin{equation}
N_{2}(F;B)\ll_{\ep}B^{3/\sqrt{d}+2/(d-1)+\ep}\log^{5}||F||.
\end{equation}
To bound $N_{1}(F;B)$, we observe, by Theorem 12, that there are 
$O(1)$ curves $C$ remaining.  Lines defined over $\Q$ are excluded, by
definition of $N_{1}(F;B)$, and other lines contribute $O(1)$ each.  
Thus Theorem 5, with
$2\leq\delta\leq d-2$, provides a bound $O_{\ep}(B^{1+\ep})$ for
each of the remaining curves, whence
\begin{equation}
N_{1}(F;B)\ll_{\ep}B^{1+\ep}+B^{3/\sqrt{d}+2/(d-1)+\ep}\log^{5}||F||.
\end{equation}
Finally, we note that the analogues in $\Proj^{3}$ of the 
bounds (1.3) and (1.4) show that curves of genus
at least $1$, and degree at most $d-2$, contribute
$O_{\ep,F}(B^{\ep})$, so that (8.1) implies (1.16).  As in Corollary 1,
we may in fact restrict attention to curves defined over the rationals
in applying the analogues of (1.3) and (1.4).

As in Section 6 we have to eliminate the factor $\log^{5}||F||$ from
(8.1) and (8.2), and we apply the same technique.  In the case of (8.2) the
argument is exactly as before.  For the bound (8.1)
we take $\x^{(1)},\ldots,\x^{(N)}$ to be the
complete set of solutions of $F(\x)=0$ in the region $|\x^{(i)}|\ll
B$, excepting any that lie on curves of degree at most $d-2$ in the surface.
This time we have must have $\delta\geq d-1$ for the curves $C$
arising in the second case, so that
$$N=N_{2}(F;B)\ll_{\ep}B^{2/(d-1)+\ep}.$$ 
As before this is sufficient.

\section{Sums of $3$ powers}

This section will be devoted to the proof of Theorem 13.  We shall take
$F(\x)=x_{1}^{d}+x_{2}^{d}+x_{3}^{d}-Nx_{4}^{d}$, and consider points
with $0<x_{1},x_{2},x_{3}\leq B$ and $x_{4}=1$.  Such points
have $\x\in Z_{4}$, and lie in
the box $|x_{i}|\leq B_{i}$, with $B_{1}=B_{2}=B_{3}=B$ and
$B_{4}=1$, so that we have $V=B^{3}$ and $T=B^{d}$, in the notation of
Theorem 14.  An application of Theorem 14 therefore shows that our
points lie on one of $O_{\ep}(B^{2/\sqrt{d}+\ep})$ curves, each
having degree $O_{\ep}(1)$.  If such a curve has degree
$D\geq d-1$, Theorem 5 then shows that there are
$O_{\ep}(B^{2/D+\ep})$ corresponding points.  This produces a total
of $O_{\ep}(B^{\theta+\ep})$ points, with
$\theta=2/\sqrt{d}+2/(d-1)$, as in Theorem 13.  This is acceptable.

We now turn to curves $C$ of degree at most $d-2$.  Let
$\theta:\Proj^{3}\rightarrow\Proj^{3}$ be the map
$$\theta(x_{1},x_{2},x_{3},x_{4})=(x_{1},x_{2},x_{3},N^{1/d}x_{4}).$$
Then $\theta(C)$ is a curve of degree at most $d-2$, lying in
the nonsingular surface $S$, given by
$y_{1}^{d}+y_{2}^{d}+y_{3}^{d}-y_{4}^{d}=0$.
According to Theorem 14 there are $O(1)$ such curves,
$C_{1},\ldots,C_{t}$, say.  Clearly $t$ and the curves $C_{i}$ depend
only on $d$ and not on $N$.  Since the point $(0,0,0,1)$ does not lie
on the surface $x_{1}^{d}+x_{2}^{d}+x_{3}^{d}-Nx_{4}^{d}=0$, it cannot
lie on the curve $C$, so that the projection 
$$\pi(x_{1},x_{2},x_{3},x_{4})=(x_{1},x_{2},x_{3})$$
is a regular map from $C$ to a curve $\pi(C)$ in $\Proj^{2}$.
Similarly $\pi$ is a regular map from each $C_{i}$ to a curve
$\pi(C_{i})=\Gamma_{i}$, say in $\Proj^{2}$.  It is
clear from the definitions that $\pi\theta=\pi$.  Thus if
$\theta(C)=C_{i}$, then $\pi(C)=\pi(C_{i})=\Gamma_{i}$.  It follows that if
$(x_{1},x_{2},x_{3},1)$ lies on a curve $C$, then
$(x_{1},x_{2},x_{3})$ lies on one of $O(1)$ curves $\Gamma_{i}$,
which are independent of $N$.  If $\Gamma_{i}$ is not defined over
$\Q$ then it has $O(1)$ rational points, by Corollary 1.  Similarly if
the genus of $\Gamma_{i}$ is $2$ or
more, we deduce from Faltings' theorem (1.4) that there are $O(1)$
rational points.  In these cases $(x_{1},x_{2},x_{3})$
is a scalar multiple of one of $O(1)$ points.  At most one such
scalar multiple can satisfy the additional relation
$x_{1}^{d}+x_{2}^{d}+x_{3}^{d}=N$. If $\Gamma_{i}$ has genus $1$, it 
follows from (1.3) that
there are $O_{\ep}(B^{\ep})$ possible points $(x_{1},x_{2},x_{3})$, up
to multiplication by scalars, and again there can be at most one
admissible scalar multiple for each value of $N$.  We therefore
conclude that there are $O_{\ep}(B^{\ep})$ solutions to
$x_{1}^{d}+x_{2}^{d}+x_{3}^{d}=N$, corresponding to points on curves
$C$ for which $\pi(C)=\Gamma_{i}$ has positive genus.

It remains to consider the possibility that $\pi(C)=\Gamma_{i}$ is
defined over $\Q$ and has genus zero.  We shall assume, as we clearly 
may, that the curve has infinitely many rational
points.  We proceed to show that $\Gamma_{i}$ can be parametrized.
Write the curve in affine coordinates as $f(x,y)=0$, where
$f(x,y)\in\Q[x,y]$ is absolutely irreducible.  We may clearly choose
the coordinates so that at most finitely many points lie at infinity.
According to Eichler
[5; p.\  139] there are two possibilities.  It may happen that the 
function field
$\Q(x,y)$ is a rational function field $\Q(z)$ for some $z\in\Q(x,y)$.
Alternatively, we may have $\Q(x,y)=\Q(u,v)$ where $g(u,v)=0$ 
for some quadratic
polynomial $g(u,v)\in\Q[u,v]$ having the property that $g(a,b)=0$ has
no rational solutions $a,b$.  In this second case we may write $u$ and
$v$ as rational functions 
$$u(x,y)=U(x,y)/W(x,y),\;\;\;v(x,y)=V(x,y)/W(x,y),$$
where $U(x,y),V(x,y),W(x,y)\in\Z[x,y]$ and $f(x,y)\nmid W(x,y)$.  Now,
if we have a rational point $(a_{1},a_{2},a_{3})$ on $\Gamma_{i}$ and
$a_{3}\not=0$, then $f(b_{1},b_{2})=0$ with $b_{1}=a_{1}/a_{3},\,
b_{2}=a_{2}/a_{3}$.  Thus if $W(a_{1}/a_{3},a_{2}/a_{3})\not=0$ we see
that
$$c_{1}=\frac{U(a_{1}/a_{3},a_{2}/a_{3})}{W(a_{1}/a_{3},a_{2}/a_{3})},\;\;\;
c_{2}=\frac{V(a_{1}/a_{3},a_{2}/a_{3})}{W(a_{1}/a_{3},a_{2}/a_{3})},$$
is a rational solution of $g(c_{1},c_{2})=0$.  This contradiction
would show, in this second case, that every rational point
$(a_{1},a_{2},a_{3})$ on $\Gamma_{i}$ would have to satisfy either $a_{3}=0$
or $W(a_{1}/a_{3},a_{2}/a_{3})=0$.  Since $f(x,y)\nmid W(x,y)$ this
allows only finitely many points, which again contradicts our initial
assumptions.  Thus we must be in the first case, in which 
$\Q(x,y)=\Q(z)$ for some $z\in\Q(x,y)$.

We now revert to the projective formulation of the curve $\Gamma_{i}$.
From the fact that $\Q(x,y)=\Q(z)$ we conclude that there are integral
binary forms $f_{1}(u,v),f_{2}(u,v)$ and $f_{3}(u,v)$, with no common
factor, 
such that, if $\x$ lies on $\Gamma_{i}$, then it is
proportional to $(f_{1}(u,v),f_{2}(u,v),f_{3}(u,v))$ for some
$(u,v)$.  Moreover there are coprime forms $u(\x),v(\x)
\in\Z[x_{1},x_{2},x_{3}]$ whose ratio is nonconstant on $\Gamma_{i}$,
such that the appropriate values of $u$ and $v$ may be given by 
$u=u(\x)$ and $v=v(\x)$.  Thus a
rational point $\x$ on $\pi(C)$ will be a nonzero rational scalar
multiple of $(f_{1}(u,v),f_{2}(u,v),f_{3}(u,v))$ for some primitive
$(u,v)\in\Z^{2}$, except in a finite number of
cases.  These exceptions arise when $f_{i}(u,v)=0$ for $i=1,2,3$ and 
$u=u(\x), v=v(\x)$, and hence there are $O(1)$ of them.  We
therefore see that, apart from these exceptions, the relevant points on the
curve $C$ are given by solutions of
$$\lambda^{d}(f_{1}(u,v)^{d}+f_{2}(u,v)^{d}+f_{3}(u,v)^{d})=N,\;\;\;
\lambda\in\Q,\;\;u,v\in\Z,\;\;(u,v)=1.$$
Since the
forms $f_{i}$ have no common factor there will be relations of
the type
$$\sum_{i=1}^{3}g_{i}(u,v)f_{i}(u,v)=Gu^{r},\;\;\;
\sum_{i=1}^{3}h_{i}(u,v)f_{i}(u,v)=Hv^{r}.$$
Here $g_{i}(u,v),h_{i}(u,v)$ are integral forms, and $G,H$ are nonzero integer
constants.  We note that 
$\lambda f_{i}(u,v)$ must
be integral for $i=1,2,3$, in order to produce integral values of
$x_{i}$.  Since $u$ and $v$ are coprime, it follows that the 
denominator of $\lambda$
must divide $GH$, and hence can take only $O(1)$ values.  Setting
$\lambda=\mu/\nu$ with $\mu,\nu$ coprime, we have $\mu^{d}|N$, so that
$\mu$ takes $O_{\ep}(N^{\ep})$ values.  

It remains to consider the
number of solutions in $u,v$ that the Thue equation
\begin{equation}
f_{1}(u,v)^{d}+f_{2}(u,v)^{d}+f_{3}(u,v)^{d}=\nu^{d}\mu^{-d}N
\end{equation}
may have.  We shall write $f(u,v)$ for the form on the left-hand side.
Recall that $\mu,\nu$ and the form $f$
may be considered as fixed.  Suppose firstly that $f$ has two 
distinct rational factors, $f'$ and $f''$, say, both irreducible over
$\Q$.  Then
$f'(u,v)=N'$ and $f''(u,v)=N''$ for certain factors $N',N''$ of 
$\nu^{d}\mu^{-d}N$.  These two equations determine $O(1)$ values of
$u,v$, by elimination, so that (9.1) has $O_{\ep}(N^{\ep})$ solutions.
There remains the possibility that $f$ is a power of an irreducible
form $f'$ say, in which case we have to consider solutions of an
equation $f'(u,v)=N'$.  If $f'$ has degree $3$ or more we can apply the result
of Lewis and Mahler [25], which shows that there are
$O(A^{\omega(N')})$ such solutions, with a constant $A$ depending
only on $f'$.  Our construction shows that this latter form is one
of a finite set, independent of $N$.  There are therefore
$O_{\ep}(N^{\ep})$ solutions in this case.  

When $f'$ has degree two
the equation $f'(u,v)=N'$ will have $O_{\ep}(N^{\ep})$ solutions,
providing that
the variables $u,v$ can be bounded by powers of $N$.  Since we assumed
that the forms $f_{i}$ had no common factor, we may take $f_{1}$, say,
to be coprime to $f'$.  However if $f'(u,v)\ll N$ then $u/v-\alpha\ll
N/|v|$, for some root $\alpha$ of $f'(X,1)$ (unless $v=0$).  Similarly, from
$f_{1}(u,v)\ll N$ we have $u/v-\beta\ll N/|v|$ for some root $\beta$ of
$f_{1}(X,1)$, unless $v=0$. Since $f'$ and $f_{1}$ are coprime we will
have $\alpha\not=\beta$, and hence $N/|v|\gg 1$.  It follows that
$v\ll N$, whether or not $v\not=0$.  Similarly we have $u\ll N$.  This
gives us the necessary bounds on $u$ and $v$.

We have therefore shown that the equation (9.1)
has $O_{\ep}(N^{\ep})$ solutions, except possibly when the form
$f(u,v)$ on the left-hand side is a constant multiple of a power of
a rational
linear function $L(u,v)$, say.  In this last case, we may make an
appropriate linear change of variable, invertible over $\Z$, so that we
actually have $f(u,v)=cv^{dk}$ for some $k\in\N$ and some
nonzero integer constant $c$.
We have therefore to ask whether an identity of the form
\begin{equation}
f_{1}(u,v)^{d}+f_{2}(u,v)^{d}+f_{3}(u,v)^{d}=cv^{dk}
\end{equation}
is possible, with coprime integral forms $f_{i}$ of degree $k$.  If 
$f_{i}(u,v)$ has leading term $a_{i}u^{k}$ in $u$, we will have $a_{1}^{d}+
a_{2}^{d}+a_{3}^{d}=0$.  Thus $d$ must be odd, and in view of Wiles'
proof of Fermat's last theorem [31], we may assume that $a_{2}=-a_{1}$
and $a_{3}=0$.  Since we then have $v^{d}|f_{3}(u,v)^{d}$ and
$v^{d}|cv^{dk}$ we conclude that
$v^{d}|f_{1}(u,v)^{d}+f_{2}(u,v)^{d}$.  However
\begin{eqnarray}
&&f_{1}^{d-1}-f_{1}^{d-2}f_{2}+\ldots-f_{1}
f_{2}^{d-2}+f_{2}^{d-1}\\
&&\qquad\quad \equiv \lrp{a_{1}^{d-1}-a_{1}^{d-2}a_{2}+\ldots-a_{1}a_{2}^{d-2}+a_{2}^{d-1}}
u^{(d-1)k}\mod{v},\nonumber
\end{eqnarray}
and $a_{1}^{d-1}-a_{1}^{d-2}a_{2}+\ldots-a_{1}a_{2}^{d-2}+a_{2}^{d-1}
=da_{1}^{d-1}$,
since $a_{2}=-a_{1}$.  Moreover $a_{1}\not=0$, for otherwise the forms
$f_{i}$ would not be coprime. It follows that the expression (9.3) is
coprime to $v$, and hence that $v^{d}|
f_{1}(u,v)+f_{2}(u,v)$.  We cannot have $f_{2}(u,v)=-f_{1}(u,v)$ since
the parametrization could not then produce solutions in positive
integers.  It follows that the degree $k$ of the forms $f_{i}$ must be
at least $d$.  The relations (9.1) and (9.2) show that $\mu,\nu$ and $N$ 
determine
$|v|$, and that $v\ll N^{1/dk}\ll B^{1/k}\ll B^{1/d}$.  Then, from
$f_{1}(u,v)\ll N^{1/d}$, we deduce that $u-\alpha v\ll N^{1/dk}$, for
some factor $u-\alpha v$ of $f_{1}(u,v)$.  Thus $u\ll |\alpha v|+N^{1/dk}\ll 
B^{1/d}$.  It therefore 
follows finally that
a curve $\Gamma_{i}$ of genus $0$ contributes $O(B^{1/d})$ points, which
is satisfactory for Theorem 13.

\vglue12pt
{\vbox{\baselineskip10pt
\eightsc  Mathematical Institute, Oxford University, Oxford, UK

{\eightpoint {\it E-mail address\/}: rhb@maths.ox.ac.uk}}}
 
\vglue-36pt


\vfill
\pagebreak

 \centerline{\bf Appendix}
\vglue12pt
\centerline{ By {\tensc J.-L. Colliot-Th\'el\`ene}}
\vglue12pt

\def\titleheadline#1{\def\one{#1}\ifx\one\empty\else
\gdef\thetitle{{\frenchspacing%
\let\\ \relax
{#1}}}\fi}
\newif\ifshort
\def\shortname#1{\global\shorttrue\xdef
\theauthors{{\eightsc\uppercase{J.-L. Colliot-Th\'el\`ene}}}}
\let\shorttitle\titleheadline
\shorttitle{Rational points on curves and surfaces}

Throughout this appendix the ground field is algebraically
closed of characteristic zero.
 
\specialnumber{1}\proclaim{Proposition}
Let $X \subset \Proj^{3}$ be a smooth projective surface{\rm .}
If the number of reduced and irreducible curves $C$ of degree $d$ lying on
$X$ is not finite{\rm ,} then there exists such a curve $C \subset X$
whose self\/{\rm -}\/intersection $(C.C)$ is nonnegative{\rm .}
\endproclaim

For the proof we let $G$ denote the open
set of the Hilbert scheme of curves  in $\Proj^{3}$
corresponding to integral (i.e. reduced and irreducible) curves
of degree $d$.
This is a scheme of finite type over the ground field.
Let $W \subset G \times \Proj^{3}$ be the
reduced closed subset whose points are pairs $(c,x)$
where $c$ is a point with associated curve $C$
and $x \in \Proj^{3}$ lies on $C$. For each $c \in G$,
the fibre of $W \to G$ above $c$ is just the
projective integral curve $C$.

Let $Z \subset G \times X$ be the trace of $W$
on $G \times X$.
The projection map  $Z \to G$ is projective.
Its image is a closed subset $F \subset G$.
A curve $C$ with associated point $c$ lies on $X$ if and
only if $c$ belongs to $F$, and the inverse image of
such a $c$ in $Z$ is precisely the curve $C \subset X$.
If each component of
$F$ is of dimension zero, then $F$ is finite, and there are
only finitely many curves $C$ of degree $d$ lying on
$X$. If that is not the case,
then $F$ contains at least one irreducible curve $T$.
Two distinct fibres $C_0$ and $C_1$ of $Z \to G$
above rational
points of $T$ define integral  curves in
the same algebraic family on $X$. Hence $(C_0.C_0)
=(C_0.C_1) \geq 0$, as claimed.

\specialnumber{2}\proclaim{Proposition}
Let $X \subset \Proj^{3}$ be a smooth projective
surface of
degree $n$ and let $C \subset X$ be a
reduced and irreducible curve{\rm ,}  possibly singular{\rm ,}
of degree $d$ in $\Proj^{3}${\rm .}
If $n > d+1${\rm ,} then the intersection number $(C.C)$
of $C$ on the surface $X$ is strictly negative{\rm .}
\endproclaim

The canonical sheaf $K$ on $X$ is ${\cal O}_X(n-4)$.
The formula for the arithmetic genus of $C \subset X$
is well known to be
$$2p_a(C)-2=(C.C)+(C.K)=(C.C)+(n-4)d
$$
(see [3; Chapter V, Exercise 1.3, p.\ 366]).
Here one should recall that $p_a(C)$ is by definition
the dimension of $H^1(C,{\cal O}_C)$.

Suppose first that $C$ is contained in a plane $\Proj^2 \subset \Proj^3$.
Then the formula for the arithmetic genus of $C \subset \Proj^2$
is
$$2p_a(C)-2=d(d-3).$$
\noindent Thus
$$ (C.C)+(n-4)d = d(d-3),$$
and hence
$$(C.C)= d(d+1-n)$$

\noindent 
is strictly negative as soon as $n > d+1$.

Suppose that $C$ is not contained in a plane, hence in
particular $d \geq 3$. In classical
parlance, such a curve is called nondegenerate. For such
curves we have the Castelnuovo bound:
If $d$ is even,   $p_a(C) \leq (d^2/4) - d +1$; if $d$ is odd,
  $ p_a(C) \leq (d^2-1)/4 - d +1$.

Comparing with the formula for $p_a(C)$, we find
that  $(C.C)<0$ if\break $n > d/2 + 2$, whether $d$ is even or odd.  This
completes the proof of the proposition.

Let us comment on the Castelnuovo bound.
Standard textbooks give the
Calstelnuovo bound for the genus $g(C)$ of
smooth nondegenerate curves $C \subset \Proj^3$, see  [3;
Theorem IV.6.4, p.\  351], for example. (See also [2; p.\  116].)
However the whole argument is valid for reduced, irreducible, local
complete intersection curves such as the ones under consideration here.
Indeed, the proof uses the Riemann-Roch theorem for the curve $C$,
and it uses a General Position Lemma, to
the effect that the section of $C$ by a sufficiently general
plane in $\Proj^3$ consists of $d$ distinct points, no three of
which are on a line.
  The General Position Lemma is valid for
singular (reduced, irreducible) curves ([2; p.\  109]; see also [4]
and references therein -- this reference was pointed out to me by D.
Perrin).
A proof of the Riemann-Roch theorem for (possibly singular) curves
given as divisors on a surface (such curves are automatically local
complete intersections)  may be found in [1; Chapter VII, Section 1]
(combine theorem (1.4), theorem (1.15) and Remark (1.17)).

\specialnumber{3}\proclaim{Proposition}
If for each nonsingular surface $X \subset \Proj^3$ of degree $n$ the
number of curves $C$  of degree $d$ lying on $X$ is finite{\rm ,} then there is
an integer
$N(n,d)$ such that any nonsingular surface $X \subset \Proj^3$ of degree
$n$ contains at most $N(n,d)$ curves of degree $d${\rm .}
\endproclaim

Let $W \subset \Proj(H^0(\Proj^3,{\cal O}(n)))$ be the open set
corresponding to nonsingular surfaces of degree $n$. Let
$Z \subset G \times W$ be the closed set whose points
are pairs of points $(c,f)$ with $c$ corresponding to
an integral curve $C$ of degree $d$ lying on the surface
$X$ defined by $f$. The projection map $Z \to W$ is
a proper morphism. By hypothesis, any fibre of this
morphism is finite, i.e. the morphism $Z \to W$
is quasi-finite.  Being both proper and quasi-finite,
the morphism $Z \to W$  is finite, see  [3; Exercise
III.11.2, p.\  280].
This implies the existence of
an integer $N$ such that any fibre has at most
$N$ points.

\demo{{R}emark}  Computing dimensions, one sees that for $n$ big
enough with respect to the degree, the general surface of degree $n$
contains no (reduced, irreducible) curve of degree $d$ at all.
\enddemo

 \pagebreak
 
Gathering the three propositions together, we conclude as follows.

\nonumproclaim{Theorem}
For each pair $n,d$ of positive
integers with $n > d+1${\rm ,} there exists an integer
$N(n,d)$ such that for any smooth projective
surface $X \subset \Proj^3$ of degree $n${\rm ,} there are
at most $N(n,d)$ reduced and irreducible curves of
degree $d$ lying on $X${\rm .}
\endproclaim

\demo{{R}emark}  From this one may conclude an analogous
result where one omits the condition `reduced
and irreducible'. Indeed an effective Cartier divisor
  $C \subset X$
of degree $d$ in $\Proj^3$
defines a divisor $\sum_in_iC_i$ with $n_i>0$, where
each $C_i$ is an integral curve of degree $d_i$, and the sum
$\sum_in_id_i$ is equal to $d$.
\enddemo
{\vbox{\baselineskip10pt
\eightsc  Universit\'e de Paris-Sud, Orsay, France 

{\eightpoint {\it E-mail address\/}:  colliot@math.u-psud.fr}}}
 
\vglue-24pt

\centerline{\ninerm (Received November 17, 2000)}

\end{document}

\end{document}